\documentclass[12pt]{article}
\usepackage{amsmath}
\usepackage{graphicx,psfrag,epsf}
\usepackage{enumerate}
\usepackage{natbib}
\usepackage{url} 
\usepackage{bm}
\usepackage{algorithm}
\usepackage{subfig}
\usepackage{amsfonts}
\usepackage{blkarray}
\usepackage{url}

\graphicspath{{./}{./figs/}}

  \def\clap#1{\hbox to 0pt{\hss#1\hss}}

\providecommand{\mat}[1]{\bm{#1}}%
\renewcommand{\vec}[1]{\mathbf{#1}}

\providecommand{\mA}{\ensuremath{\mat{A}}}

\providecommand{\mC}{\ensuremath{\mat{C}}}
\providecommand{\mD}{\ensuremath{\mat{D}}}

\providecommand{\mU}{\ensuremath{\mat{U}}}

\providecommand{\mW}{\ensuremath{\mat{W}}}

\providecommand{\vu}{\ensuremath{\vec{u}}}
\providecommand{\vv}{\ensuremath{\vec{v}}}

\providecommand{\vx}{\ensuremath{\vec{x}}}
\providecommand{\vy}{\ensuremath{\vec{y}}}
\providecommand{\vz}{\ensuremath{\vec{z}}}

\newcommand{\hLambda}{\hat{\Lambda}}

\newcommand{\hmW}{\hat{\mW}}

\newcommand{\bPi}{\bm{\Pi}}

\newcommand{\bmat}[1]{\begin{bmatrix}#1\end{bmatrix}}

\newcommand{\rank}{\operatorname{rank}}

\newcommand{\unit}[1]{\text{#1}}

\newcommand{\uavg}{u_{\text{avg}}}
\newcommand{\Bind}{B_{\text{ind}}}
\renewcommand{\Re}{Re}
\newcommand{\Ha}{Ha}

\newcommand{\tL}{\text{L}}
\newcommand{\tT}{\text{T}}
\newcommand{\tM}{\text{M}}
\newcommand{\tC}{\text{C}}
\newcommand{\matindex}[1]{\mbox{\scriptsize#1}}

\begin{document}

\def\spacingset#1{\renewcommand{\baselinestretch}%
{#1}\small\normalsize} \spacingset{1}


\title{\bf Dimension reduction in MHD power generation models: dimensional analysis and active subspaces}
\author{
Andrew Glaws and 
Paul G.~Constantine \\
Department of Applied Mathematics and Statistics \\
Colorado School of Mines, Golden, CO
\vspace{.2cm} \\
John Shadid and Timothy M.~Wildey \\
Sandia National Laboratories, Albuquerque, NM}

\maketitle

\bigskip
\begin{abstract}
Magnetohydrodynamics (MHD)---the study of electrically conducting fluids---can be harnessed to produce efficient, low-emissions power generation. Today, computational modeling assists engineers in studying candidate designs for such generators. However, these models are computationally expensive, so studying the effects of the model's many input parameters on output predictions is typically infeasible. We study two approaches for reducing the input dimension of the models: (i) classical dimensional analysis based on the inputs' units and (ii) active subspaces, which reveal low-dimensional subspaces in the space of inputs that affect the outputs the most. We also review the mathematical connection between the two approaches that leads to consistent application. The dimension reduction yields insights into the driving factors in the MHD power generation models. We study both the simplified Hartmann problem, which admits closed form expressions for the quantities of interest, and a large-scale computational model with adjoint capabilities that enable the derivative computations needed to estimate the active subspaces. 
\end{abstract}

\noindent%
\vfill
\hfill {\tiny technometrics tex template (do not remove)}

\newpage
\spacingset{1.45} 

\section{Introduction}
\label{sec:intro}

Magnetohydrodynamics (MHD) is an area of physics concerned with electrically conducting fluids~\citep{Cowl:57}; its various models are found in disparate fields from geophysics to fusion energy. The US Department of Energy\textsc{\char13}s investment in MHD for power generation with low emissions generators dates back to the 1960s. Interest waned in the mid 1990s after a proof-of-concept program revealed several technical challenges to scaling and integrating MHD-based components into a practical power generator~\citep{gao93}. However, in the last two decades, several of these challenges have been addressed via other investments; a recent workshop at the National Energy Technology Laboratory brought together current MHD researchers from labs, academia, and industry to assess the state of the art in MHD power generation. In particular, given two decades of supercomputing advances, the workshop report calls for improved simulation tools that exploit modern supercomputers to build better MHD models and integrate them into full scale generator designs~\citep{mhdworkshop}.

Our recent work has developed scalable simulations for resistive MHD models~\citep{Lin2010,Shadid2010,Shadid2016,Sondak}. To incorporate such simulations into generator design, a designer must understand the sensitivities of model predictions to changes in model input parameters. With this in mind, we have developed adjoint capabilities in the MHD simulation codes that enable computation of derivatives of output quantities of interest with respect to input parameters~\citep{Shadid2016a}. 

Derivative-based sensitivity analysis can reveal a low-dimensional parameterization of the map from model inputs to model predictions by identifying combinations of parameters whose perturbations change predictions the most; less important combinations can be safely ignored when assessing uncertainty in predictions or employing the model for design. Such dimension reduction may enable exhaustive parameter studies not otherwise feasible with the expensive simulation model. In this paper, we exploit the adjoint-based gradient capabilities to uncover an \emph{active subspace}~\citep{asmbook} in each of two quantities of interest (the average flow velocity and the induced magnetic field) from a laminar flow MHD model as a proof-of-principle demonstration. An active subspace is the span of a set of directions in the model's input space; perturbing the inputs along these directions changes the quantity of interest more, on average, than perturbing the inputs along orthogonal directions. 

Our recent work has revealed active subspaces in physics-based simulations for solar cell models~\citep{constantine2015discovering}, integrated hydrologic models~\citep{Jefferson2015}, and multiphysics scramjet models~\citep{Constantine2015}. We distinguish the current application to MHD by connecting the active subspace-based dimension reduction to dimensional analysis---a standard tool in science and engineering for reducing the number of parameters in a physical model by examining the physical quantities' units. The celebrated Buckingham Pi Theorem (see~\citet[Chapter 1]{Barenblatt1996}) loosely states that a physical system with $m$ parameters whose derived units are products of powers of $k<m$ fundamental units (e.g., meters, seconds, etc.) can be written in terms of $m-k$ unitless quantities. Dimensional analysis can reveal fundamental insights into the physical model before running any simulations or conducting any experiments---merely by examining units. Our recent work has revealed a mathematical connection between standard dimensional analysis and dimension reduction with active subspaces~\citep{constantine2016many}; namely, the dimensional analysis provides an upper bound on the number of important input space directions without any computation. We exploit this connection to gain insight into the driving factors of two MHD models. 

The remainder of the paper proceeds as follows. In Section \ref{sec:dimred} we review dimensional analysis, active subspaces, and the connections between dimensional analysis and active subspaces. In Section \ref{sec:mhddr} we apply dimensional analysis to the governing equations of MHD to determine the number of dimensionless quantities that affect the system. In Section \ref{sec:mhdas}, we estimate active subspaces in two MHD models: (i) the Hartmann problem, which is a simplified description of an MHD duct flow that has direct relevance to MHD generators and (ii) large-scale steady state laminar flow simulation of an MHD generator in three spatial dimensions. In both cases, we comment on the insights revealed by comparing dimensional analysis to the gradient-based, computationally discovered active subspace.

\section{Dimension reduction methods}
\label{sec:dimred}

Reducing a system's dimension can both yield insight into its fundamental characteristics and enable parameter studies that are otherwise infeasible when the dimension is too large. The type of dimension reduction in the present study applies to systems of the form 
\begin{equation}
\label{eq:fx}
y \;=\; f(\vx), \qquad y\in\mathbb{R},\qquad \vx\in\mathbb{R}^m,
\end{equation}
where the function $f$ maps $\mathbb{R}^m$ to $\mathbb{R}$. The scalar $y$ is the system's quantity of interest, and $\vx$ is a vector of input parameters. Loosely speaking, we seek a parameterization of \eqref{eq:fx} with fewer than $m$ parameters. 

We assume the components of $\vx$ are independent, which distinguishes the type of dimension reduction methods we consider. In particular, there is no covariance in the components of $\vx$, so covariance based dimension reduction---such as principal component analysis---is not appropriate. Instead, we seek to reduce the dimension in a way that takes advantage of structure in $f(\vx)$ from \eqref{eq:fx}. One example of such dimension reduction is \emph{sufficient dimension reduction} in statistical regression~\citep{cook2009regression}, where one seeks a subspace in the space of $m$ predictors such that the response's statistics are unchanged. In contrast to the regression problem, the system \eqref{eq:fx} is deterministic; in other words; we do not assume that $y$ is corrupted by random noise. 

We consider two types of dimension reduction for \eqref{eq:fx}. The first is \emph{dimensional analysis}, which is a mature tool for reducing the number of variables in a physical system by examining the quantities' units. The second is based on an \emph{active subspace} in $f(\vx)$, which is a subspace of $\mathbb{R}^m$ constructed such that perturbing $\vx$ along the active subspace changes $y$ more, on average, than perturbing $\vx$ orthogonally to the active subspace. The former is analytical while the latter is computational. In other words, dimensional analysis follows from the quantities' units and can typically be performed without the aid of a computer in small systems. In contrast, the active subspace is estimated by computing $y$ and its gradient at several values of $\vx$. 

\subsection{Dimensional analysis and Buckingham Pi Theorem}

To apply dimensional analysis, we refer to \eqref{eq:fx} as a \emph{physical law}, which implies that the inputs and output are accompanied by physical units. These units are derived from $k \leq m$ base units, denoted generically as $L_1, \dots ,L_k$, which contain a subset of the seven SI units; Table \ref{tab:units} shows the SI base units.

\begin{table}[ht]
\centering
\caption{International system of units (SI). Table taken from~\citet{nistunits}.}
\begin{tabular}{lll}
Base quantity & Name & Symbol\\
\hline
length & meter & $\unit{m}$\\
mass & kilogram & $\unit{kg}$\\
time & second & $\unit{s}$\\
electric current & ampere & $\unit{A}$\\
thermodynamic temperature & kelvin & $\unit{K}$\\
amount of a substance & mole & $\unit{mol}$\\
luminous intensity & candela & $\unit{cd}$
\end{tabular}
\label{tab:units}
\end{table}

The \emph{unit function} of a quantity, denoted by square brackets, returns the units of its argument. For example, if $y$ is a velocity, then $[y]=\unit{m}\cdot\unit{s}^{-1}$. If a quantity is unitless, then its unit function returns 1. We can generically write the units of $\vx=[x_1,\dots,x_m]^T$ and $y$ as
\begin{equation} 
\label{eq:units}
\left[ x_j \right] = \prod_{i=1}^k L_i^{d_{i,j}} , \qquad 
\left[ y \right] = \prod_{i=1}^k L_i^{u_{i}}.
\end{equation}
In words, the units of each physical quantity can be written as a product of powers of the $k$ base units. We can derive unitless quantities as follows; the motivation for doing so will soon become apparent. Define the $k\times m$ matrix $\mD$ and the $k$-vector $\vu$ as
\begin{equation}
\label{eq:D}
\mD = 
\bmat{
d_{1,1} & \cdots & d_{1,m} \\ 
\vdots & \ddots & \vdots \\
d_{k,1} & \cdots & d_{k,m}
}, \qquad
\vu = \bmat{u_1 \\ \vdots \\ u_k},
\end{equation}
which contain the powers from the quantities' units in \eqref{eq:units}. We assume $\mD$ has full row rank, i.e., $\rank(\mD)=k$; a properly formulated physical model leads to $\mD$ with full row rank. Let $\vv=[v_1,\dots,v_m]^T$ satisfy $\mD\vv=\vu$, and let $\mU\in\mathbb{R}^{m\times n}$ be a basis for the null space of $\mD$, i.e., $\mD\mU=\mathbf{0}\in\mathbb{R}^{k\times n}$, where $n=m-k$. Note that the elements of $\vv$ and $\mU$ are not unique. We can construct a unitless quantity of interest $\Pi$ as
\begin{equation}
\label{eq:nondim_QoI}
\Pi \;=\; y \, \prod_{i=1}^m x_i^{-v_i},
\end{equation}
where, by constuction, $[\Pi]=1$. We similarly construct unitless parameters $\Pi_j$ as
\begin{equation} 
\label{eq:nondim_inputs}
\Pi_j \;=\; \prod_{i=1}^m x_i^{u_{i,j}}, \qquad j=1,\dots,n,
\end{equation}
where $u_{i,j}$ is the $(i,j)$ element of $\mU$. 

The Buckingham Pi Theorem~\citep[Chapter 1]{Barenblatt1996} states that a physical law \eqref{eq:fx} can be written in unitless form 
\begin{equation}
\label{eq:unitlesslaw}
\Pi \;=\; \tilde{f}(\bPi),\qquad \bPi=\bmat{\Pi_1&\cdots&\Pi_n}^T.
\end{equation}
The number of inputs in the unitless form of the physical law \eqref{eq:unitlesslaw} is $n<m$, which has reduced dimension compared to \eqref{eq:fx}. This dimension reduction may enable more accurate semi-empirical modeling of the map $\tilde{f}$ given experimental data, since a model with fewer inputs typically has fewer parameters to fit with a given data set. Additionally, the unitless quantities allow one to devise scale-invariant experiments, since scaling the units does not change the form of the unitless physical law \eqref{eq:unitlesslaw}.  

\subsection{Active subspaces}
\label{sec:as}

Assume the space of $\vx$ from \eqref{eq:fx} is equipped with a probability density function $\gamma(\vx)$. Additionally, assume that (i) $f$ is square-integrable with respect to $\gamma$ and (ii) $f$ is differentiable with square-integrable partial derivatives. Denote the gradient vector of $f$ as $\nabla f(\vx)$. Define the $m \times m$ symmetric and positive semidefinite matrix $\mC$ as
\begin{equation} 
\label{eq:C}
C \;=\; \int \nabla f(\vx)\,\nabla f(\vx)^T\, \gamma(\vx) \, d\vx.
\end{equation}
Since $\mC$ is symmetric, it admits a real eigenvalue decomposition,
\begin{equation}
\mC \;=\; \mW \Lambda \mW^T,
\end{equation}
where the columns of $\mW$ are the orthonormal eigenvectors, and $\Lambda$ is a diagonal matrix of the associated eigenvalues in descending order. 

The eigenpairs are functionals of $f$ for the given $\gamma$. Assume that $\lambda_n>\lambda_{n+1}$ for some $n<m$. Then we can partition the eigenpairs as
\begin{equation}
\Lambda = \bmat{\Lambda_1 & \\ & \Lambda_2}, \qquad
\mW = \bmat{\mW_1 & \mW_2 }, 
\end{equation}
where $\Lambda_1$ contains the first $n$ eigenvalues of $\mC$, and $\mW_1$'s columns are the corresponding eigenvectors. The \emph{active subspace} of dimension $n$ is the span of $\mW_1$'s columns; the \emph{active variables} $\vy\in\mathbb{R}^n$ are the coordinates of $\vx$ in the active subspace. The active subspace's orthogonal complement, called the \emph{inactive subspace}, is the span of $\mW_2$'s columns; its coordinates are denoted $\vz\in\mathbb{R}^{m-n}$ and called the \emph{inactive variables}.

The following property justifies the labels; see Lemma 2.2 from \cite{constantine2014active}:
\begin{equation} 
\label{eq:AS_property}
\begin{aligned}
\lambda_1 + \dots + \lambda_n 
&= \int \nabla_{\vy} f(\vx)^T  \nabla_{\vy} f(\vx)\, \gamma(\vx) \, d\vx , \\
\lambda_{n+1} + \dots + \lambda_m 
&= \int \nabla_{\vz} f(\vx)^T  \nabla_{\vz} f(\vx)\, \gamma(\vx) \, d\vx ,
\end{aligned}
\end{equation}
where $\nabla_{\vy}$ and $\nabla_{\vz}$ denote the gradient of $f$ with respect to the active and inactive variables, respectively. Since the eigenvalues are in descending order, and since $f$ is such that $\lambda_n>\lambda_{n+1}$, \eqref{eq:AS_property} says that perturbations in $\vy$ change $f$ more, on average, than perturbations in $\vz$. 

If the eigenvalues $\lambda_{n+1},\dots,\lambda_m$ associated with the inactive subspace are sufficiently small, then $f$ can be approximated by a \emph{ridge function}~\citep{pinkus2015}, which is a function that is constant along a subspace of its domain,
\begin{equation} 
\label{eq:AS_Red}
f(\vx) \;\approx\; g(\mW_1^T \vx),
\end{equation}
where $g$ maps $\mathbb{R}^n$ to $\mathbb{R}$. If we are given function evaluations $(\vx_i,f(\vx_i))$, then constructing $g$ may be feasible when constructing an approximation in all $m$ variables may not be---since $g$ is a function of $n<m$ variables. Thus, the active subspace allows dimension reduction in the space of $\vx$.

In practice, we estimate $\mC$ from \eqref{eq:C} with numerical integration; high-accuracy Gauss quadrature may be appropriate when the dimension of $\vx$ is small and the integrands are sufficiently smooth. \cite{constantine2015computing} analyzes the following Monte Carlo method. Draw $M$ independent samples $\vx_i$ according to the given $\gamma(\vx)$ and compute
\begin{equation}
\mC \;\approx\; 
\hat{\mC} \;=\; 
\frac{1}{M} \sum_{i=1}^M \nabla f(\vx_i)\, \nabla f(\vx_i)^T
\;=\;
\hmW\hLambda\hmW^T.
\end{equation}
The eigenpairs $\hLambda$, $\hmW$ of $\hat{\mC}$ estimate those of $\mC$. \cite{constantine2015computing} analyzes how large $M$ must be to ensure quality estimates; the analysis supports a heuristic of $M=\alpha(\delta)\,k\,\log(m)$ samples to estimate the first $k$ eigenvectors within a relative error $\delta$, where $\alpha(\delta)$ is typically between 2 and 10. Let $\varepsilon$ denote the distance between the true subspace and its Monte Carlo estimate defined as 
\begin{equation}
\label{eq:subdist}
\varepsilon \;=\; \|\mW_1\mW_1^T - \hmW_1\hmW_1^T\|,
\end{equation}
where $\hmW_1$ contains the first $n$ columns of $\hmW$. Corollary 3.7 from~\cite{constantine2015computing} shows that 
\begin{equation}
\label{eq:gap}
\varepsilon \;\leq\; \frac{4\,\lambda_1\,\delta}{\lambda_n - \lambda_{n+1}}.
\end{equation}
Equation \eqref{eq:gap} shows that a large eigenvalue gap $\lambda_n-\lambda_{n+1}$ implies that the active subspace can be accurately estimated with Monte Carlo. Therefore, the practical heuristic is to choose the dimension $n$ of the active subspace according to the largest gap in the eigenvalues. Additionally, \cite{constantine2014active} analyzes the effect of using estimated eigenvectors in the ridge function approximation \eqref{eq:AS_Red}. 

\subsection{Connecting dimensional analysis to active subspaces}
\label{sec:connections}

Next we study the relationship between these two dimension reduction techniques; this closely follows the development in \cite{constantine2016many}. Dimensional analysis produces the unitless physical law \eqref{eq:unitlesslaw} by defining unitless quantities \eqref{eq:nondim_QoI} and \eqref{eq:nondim_inputs} as products of powers of the dimensional quantities. On the other hand, the coordinates of the active subspace can be written as linear combinations of the system's inputs, i.e., $\vy=\mW_1^T\vx$ and $\vz=\mW_2^T\vx$ for $\vx\in\mathbb{R}^m$. Thus, these two techniques are connected via a logarithmic transformation of the input space. Combining \eqref{eq:nondim_QoI}, \eqref{eq:nondim_inputs}, and \eqref{eq:unitlesslaw},
\begin{equation} 
\label{eq:connection}
\begin{aligned}
y \prod_{i=1}^m x_i^{-v_i}
&=
\tilde{f} \left( \prod_{i=1}^m x_i^{u_{i,1}} , \dots , \prod_{i=1}^m x_i^{u_{i,n}} \right) \\
&= 
\tilde{f} \left( \exp \left( \log \left( \prod_{i=1}^m x_i^{u_{i,1}} \right) \right) , \dots , \exp \left( \log \left( \prod_{i=1}^m x_i^{u_{i,n}} \right) \right) \right) \\
&= 
\tilde{f} \left( \exp \left( \sum_{i=1}^m u_{i,1} \log(x_i) \right) , \dots , \exp \left( \sum_{i=1}^m u_{i,n} \log(x_i) \right) \right) \\
&= 
\tilde{f} \left( \exp \left( \vu_1^T \log(\vx) \right) , \dots , \exp \left( \vu_n^T \log(\vx) \right) \right)
\end{aligned}
\end{equation}
where $\log(\vx)$ returns an $m$-vector with the log of each component. Then we can rewrite $y$ as
\begin{equation}
\label{eq:dimridge}
\begin{aligned}
y
&=
\exp\left(\vv^T\log(\vx)\right)\cdot\tilde{f}\left( \exp \left( \vu_1^T \log(\vx) \right) , \dots , \exp \left( \vu_n^T \log(\vx) \right) \right)\\
&=
\tilde{g}(\mA^T\log(\vx)),
\end{aligned}
\end{equation}
where 
\begin{equation}
\mA \;=\; \bmat{\vv & \vu_1 & \cdots & \vu_n} 
\;\in\; \mathbb{R}^{m\times(n+1)},
\end{equation}
and $g:\mathbb{R}^{n+1}\rightarrow\mathbb{R}$. In other words, the unitless physical law can be transformed into a ridge function of the logs of the physical inputs. Compare \eqref{eq:dimridge} to the form of the ridge approximation \eqref{eq:AS_Red}. 

Since the physical law can be written as a ridge function, its active subspace is related to the coefficient matrix $\mA$. Let $\tilde{\vx}=\log(\vx)$. By the chain rule,
\begin{equation}
\nabla_{\tilde{\vx}} \tilde{g}(\mA^T\tilde{\vx})
\;=\;
\mA\,\nabla \tilde{g}(\mA^T\tilde{\vx}),
\end{equation}
where $\nabla_{\tilde{\vx}}$ denotes the gradient with respect to $\tilde{\vx}$, and $\nabla \tilde{g}$ is the gradient of $\tilde{g}$ with respect to its arguments. Assume $\tilde{\gamma}(\tilde{\vx})$ is a probability density function on the space of $\tilde{\vx}$. Then,
\begin{equation}
\label{eq:Ctilde}
\int \nabla_{\tilde{\vx}} \tilde{g}(\mA^T\tilde{\vx})\,\nabla_{\tilde{\vx}} \tilde{g}(\mA^T\tilde{\vx})^T \,\tilde{\gamma}(\tilde{\vx})\,d\tilde{\vx}
\;=\;
\mA\,\left(
\int \nabla \tilde{g}(\mA^T\tilde{\vx})\,\nabla \tilde{g}(\mA^T\tilde{\vx})^T \,\tilde{\gamma}(\tilde{\vx})\,d\tilde{\vx}
\right)\,\mA^T.
\end{equation}
The first column of $\mA$ is not in the null space of $\mD$ from \eqref{eq:D}, and its remaining columns are a basis for $\mD$'s null space. Therefore, $\mA$ has full column rank. Then \eqref{eq:Ctilde} shows that the active subspace for the physical law---as a function of the logs of its inputs---has dimension at least $n+1$, and it is a subspace of the $\mA$'s column space. 

The connection between the dimensional analysis and the active subspace provides an upper bound on the dimension of the active subspace. However, the eigenvalues of $\mC$ from \eqref{eq:C} rank the importance of each eigenvector-defined direction. So the eigenpairs of $\mC$ reveal more about the input/output relationship than the dimensional analysis alone. 

\section{Dimensional analysis for MHD}
\label{sec:mhddr}

Next we apply dimensional analysis to the governing equations of MHD to study the number of unitless quantities affecting the system. MHD models the behavior of electrically-conducting fluids, such as ionized liquids or plasmas. The governing equations for MHD couple the Navier-Stokes equation for fluid dynamics with Maxwell's equations for electromagnetism. Under several simplifying assumptions, we can write the governing equations of steady-state MHD as
\begin{equation}
\label{MHD_Eq}
\begin{aligned}
\nabla \cdot \left[ \rho \boldsymbol{u} \otimes \boldsymbol{u} - (p_0 + p) \boldsymbol{I} - \mu \left( \nabla \boldsymbol{u} + \nabla \boldsymbol{u}^T \right) + \frac{1}{\mu_0} \left( \boldsymbol{B} \otimes \boldsymbol{B} - \frac{1}{2} || \boldsymbol{B} ||^2 \boldsymbol{I} \right) \right] = \boldsymbol{0} , \\
\nabla \cdot \left[ \boldsymbol{u} \otimes \boldsymbol{B} - \boldsymbol{B} \otimes \boldsymbol{u} - \frac{\eta}{\mu_0} \left( \nabla \boldsymbol{B} - \nabla \boldsymbol{B}^T \right) \right] = \boldsymbol{0} , \\
\nabla \cdot \boldsymbol{u} = 0 , \hspace{2em} \nabla \cdot \boldsymbol{B} = 0 ,
\end{aligned}
\end{equation}
where $\boldsymbol{u}$ is the fluid velocity, $\boldsymbol{B}$ is the magnetic field, and $p_0$ is the applied pressure.
To apply dimensional analysis, we first determine the fundamental dimensional quantities in the model. The necessary base units are $\tL = \text{length}$, $\tT = \text{time}$, $\tM = \text{mass}$, and $\tC = \text{electric current}$; see Table \ref{tab:units}. We note that $\mu_0$ is the permeability of a vacuum with units $\frac{\tM \tL}{\tT^2 \tC^2}$. This is a physical constant that does not factor into the dimensional analysis. The model's input parameters, with their units, are
\begin{itemize}
\item length, $\ell$, with $\displaystyle [ \ell ] = \tL$,
\item velocity, $v$, with $\displaystyle [ v ] = \frac{\tL}{\tT}$,
\item fluid viscosity, $\mu$, with $\displaystyle [ \mu ] = \frac{\tM}{\tL \, \tT}$,
\item fluid density, $\rho$, with $\displaystyle [ \rho ] = \frac{\tM}{\tL^3}$,
\item pressure, $p$, with $\displaystyle [ p ] = \frac{\tM}{\tL \, \tT^2}$,
\item fluid magnetic resistivity, $\eta$, with $\displaystyle [ \eta ] = \frac{\tM \, \tL^3}{\tT^3 \, \tC^2}$,
\item and magnetic field, $B$, with $\displaystyle [ B ] = \frac{\tM}{\tT^2 \, \tC}$.
\end{itemize}
The matrix $\mD$ from \eqref{eq:D} is
\begin{equation} 
\label{D_matrix}
D \;=\; \begin{blockarray}{cccccccc}
     & \matindex{$\ell$} & \matindex{$v$} & \matindex{$\mu$} & \matindex{$\rho$} & \matindex{$p$} & \matindex{$\eta$} & \matindex{$B$} \\
    \begin{block}{c[ccccccc]}
      \matindex{L} &  1 &  1 & -1 & -3 & -1 &  3 &  0 \\
      \matindex{T} &  0 & -1 & -1 &  0 & -2 & -3 & -2 \\
      \matindex{M} &  0 &  0 &  1 &  1 &  1 &  1 &  1 \\
      \matindex{C} &  0 &  0 &  0 &  0 &  0 & -2 &  -1 \\
    \end{block}
  \end{blockarray} .
\end{equation}
A basis for the null space of $\mD$ is
\begin{equation} \label{null_span}
\left\{ \begin{bmatrix} 1 \\ 1 \\ -1 \\ 1 \\ 0 \\ 0 \\ 0 \end{bmatrix} , \begin{bmatrix} 1 \\ 0 \\ -1/2 \\ 0 \\ 0 \\ -1/2 \\ 1 \end{bmatrix} , \begin{bmatrix} 0 \\ -2 \\ 0 \\ -1 \\ 1 \\ 0 \\ 0 \end{bmatrix} \right\} .
\end{equation}
Since the dimension of the null space of $D$ is 3, the Buckingham Pi Theorem states that the system depends on 3 unitless quantities; see \eqref{eq:nondim_inputs}. In this case,
\begin{equation}
\begin{aligned}
\Pi_1 
&= 
\ell^1 v^1 \mu^{-1} \rho^1 p^0 \eta^0 B^0 = \frac{\rho v \ell}{\mu}, \\
\Pi_2 
&= 
\ell^1 v^0 \mu^{-1/2} \rho^0 p^0 \eta^{-1/2} B^1 = \frac{B \ell}{\eta^{1/2} \mu^{1/2}}, \\
\Pi_3 
&= 
\ell^0 v^{-2} \mu^0 \rho^{-1} p^1 \eta^0 B^0 = \frac{p}{\rho v^2}.
\end{aligned}
\end{equation}
Therefore, a particular quantity of interest will depend on at most 4 linear combinations of the log-transformed input parameters, as shown in Section \ref{sec:connections}. 

Another way to determine the unitless quantities for the MHD governing equations \eqref{MHD_Eq} is to scale both sides appropriately. We obtain
\begin{equation}
\label{eq:nondim_MHD_Eq}
\begin{aligned}
\nabla \cdot \left[ \boldsymbol{u}^* \otimes \boldsymbol{u}^* - (p_0^* + p^*) \boldsymbol{I} - \frac{1}{Re} \left( \nabla \boldsymbol{u}^* + \nabla \boldsymbol{u}^{*^T} \right) + \frac{1}{\mu_0^*} \frac{Ha^2}{Re} \left( \boldsymbol{B}^* \otimes \boldsymbol{B}^* - \frac{1}{2} || \boldsymbol{B}^* ||^2 \boldsymbol{I} \right) \right] = \boldsymbol{0} , \\
\nabla \cdot \left[ \boldsymbol{u}^* \otimes \boldsymbol{B}^* - \boldsymbol{B}^* \otimes \boldsymbol{u}^* - \frac{1}{\mu_0^*} \left( \nabla \boldsymbol{B}^* - \nabla \boldsymbol{B}^{*^T} \right) \right] = \boldsymbol{0} , \\
\nabla \cdot \boldsymbol{u}^* = 0 , \hspace{2em} \nabla \cdot \boldsymbol{B}^* = 0 ,
\end{aligned}
\end{equation}
where the stars denote unitless versions of the physical quantities in \eqref{MHD_Eq}. Notice that \eqref{eq:nondim_MHD_Eq} depends on three unitless quantities:
\begin{itemize}
\item the Reynolds number, $\displaystyle \Re = \frac{\rho v \ell}{\mu}$,
\item the Hartmann number, $\displaystyle \Ha = \frac{B \ell}{\eta^{1/2} \mu^{1/2}}$,
\item a dimensionless pressure gradient, $\displaystyle \nabla \cdot (p_0^* \boldsymbol{I}) = \nabla p_0^* = \frac{\ell}{\rho v^2} \nabla p_0 = \frac{p}{\rho v^2}$.
\end{itemize}
Scaling the Navier-Stokes equations is consistent with the unitless quantities derived from the basis \eqref{null_span}. In other words, the unitless quantities from the Buckingham Pi analysis match those in \eqref{eq:nondim_MHD_Eq}.

\section{Active subspaces for MHD}
\label{sec:mhdas}

The dimensional analysis in the previous section, coupled with the analysis from Section \ref{sec:connections}, indicates that a scalar-valued quantity of interest from an MHD model can be written as a ridge function of four linear combinations of the log transformed inputs. Therefore, we expect numerical tests to reveal an active subspace of dimension 4 or less. We study two MHD models with active subspaces: (i) the Hartmann problem that models a simplified duct flow and (ii) a numerical model of an idealized MHD generator in three spatial dimensions. 

\subsection{Hartmann problem}
\label{sec:hart}
\begin{figure}
\begin{center}
\includegraphics[width=4in]{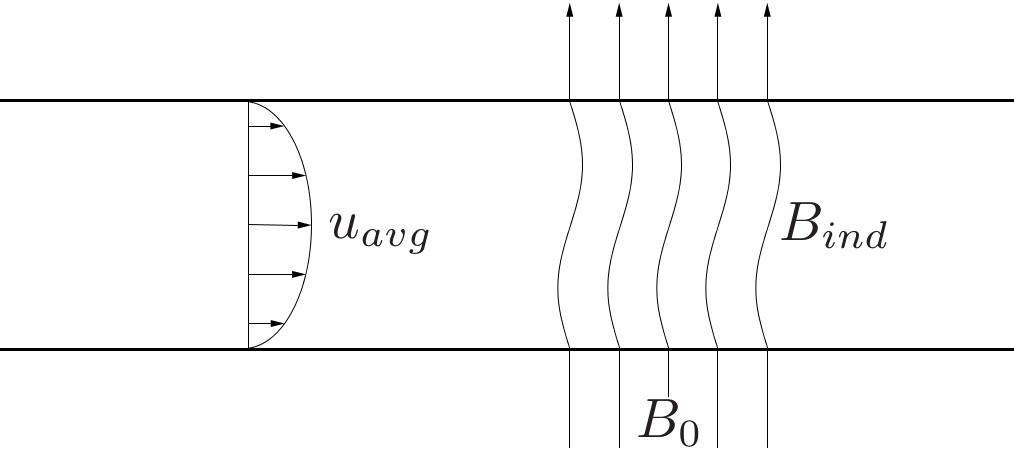}
\end{center}
\caption{Depiction of the Hartmann problem.  A magnetic fluid flows between two parallel plates in the presence of a perpendicular magnetic field.  The field acts as a damping force on the fluid while the flow induces a horizontal magnetic field. \label{fig:hartmann_prob}}
\end{figure}

The Hartmann problem is a standard problem in MHD.  It models laminar flow between two parallel plates separated by distance $2 \ell$. The fluid is assumed to be a magneto-fluid and a uniform magnetic field is applied perpendicular to the flow direction.  This magnetic field acts as a resistive force on the flow velocity while the fluid induces a magnetic field along the direction of the flow.  This can be seen in Figure \ref{fig:hartmann_prob}.

The advantage of working with the Hartmann problem is that it admits closed form analytical expressions for the quantities of interest in terms of its input parameters. This makes numerical studies of the dimension reduction much simpler. We consider two quantities of interest: (i) average flow velocity across the channel $\uavg$ and (ii) the induced magnetic field $\Bind$.  Detailed solutions to the Hartmann problem are not presented here; they may be in found in \cite{Cowl:57}.  In terms of the inputs, 
\begin{equation} 
\label{eq:Uavg_eq_dim}
\uavg \;=\; -\frac{\partial p_0}{\partial x} \frac{\eta}{B_0^2} \left( 1 - \frac{B_0 \ell}{\sqrt{\eta \mu}} \, \text{coth} \left( \frac{B_0 \ell}{\sqrt{\eta \mu}} \right) \right)
\end{equation}
and
\begin{equation} 
\label{eq:Bind_eq_dim}
\Bind \;=\; \frac{\partial p_0}{\partial x} \frac{\ell \mu_0}{2 B_0} \left( 1 - 2 \frac{\sqrt{\eta \mu}}{B_0 \ell} \, \text{tanh} \left( \frac{B_0 \ell}{2 \sqrt{\eta \mu}} \right) \right) .
\end{equation}
The expressions for the Hartmann problem quantities of interest have five input parameters: fluid viscosity $\mu$, fluid density $\rho$, applied pressure gradient $\partial p_0 / \partial x$ (where the derivative is with respect to the flow field's spatial coordinate), resistivity $\eta$, and applied magnetic field $B_0$. This differs slightly from the Buckingham Pi analysis from Section \ref{sec:mhddr}, where we considered seven inputs. Two extra quantities needed to perform dimensional analysis are fixed in the Hartmann problem. 

To connect the inputs in the Hartmann problem's quantities of interest to the notation from Section \ref{sec:as}, let
\begin{equation}
\vx \;=\; \begin{bmatrix} \log(\mu) & \log(\rho) & \log\left(\frac{\partial p_0}{\partial x}\right) & \log(\eta) & \log(B_0) \end{bmatrix}^T.
\end{equation}
To estimate the active subspace for each quantity of interest, we set $\gamma(\vx)$ from \eqref{eq:C} to be a uniform density on a five-dimensional hypercube. The ranges of each of $\vx$'s components are in Table \ref{tab:input_ints}; they are chosen to represent the expected operating conditions of an MHD generator modeled with the Hartmann problem. We estimate $\mC$ from \eqref{eq:C} with a tensor product Gauss-Legendre quadrature rule with 11 points in each dimension---a total of 161051 points. This was sufficient for 10 digits of accuracy in the eigenvalue estimates. 

\begin{table}
\centering
\caption{Indices and intervals for the parameters $\vx$ of the Hartmann problem. These intervals represent the expected operating conditions for an MHD generator modeled with the Hartmann problem.}
\begin{tabular}{llll}
Index & Name & Notation & Interval \\ \hline
1 & fluid viscosity & $\log(\mu)$ & $[\log(0.05), \log(0.2)]$ \\
2 & fluid density & $\log(\rho)$ & $[\log(1), \log(5)]$ \\
3 & applied pressure gradient & $\log\left(\frac{\partial p_0}{\partial x}\right)$ & $[\log(0.5), \log(3)]$ \\
4 & resistivity & $\log(\eta)$ & $[\log(0.5), \log(3)]$ \\
5 & applied magnetic field & $\log(B_0)$ & $[\log(0.1), \log(1)]$
\end{tabular}
\label{tab:input_ints}
\end{table}

\begin{figure}[!h]
\centering
\subfloat[Eigenvalues of $\mC$]{
\label{fig:evalsU}
\includegraphics[width=0.42\textwidth]{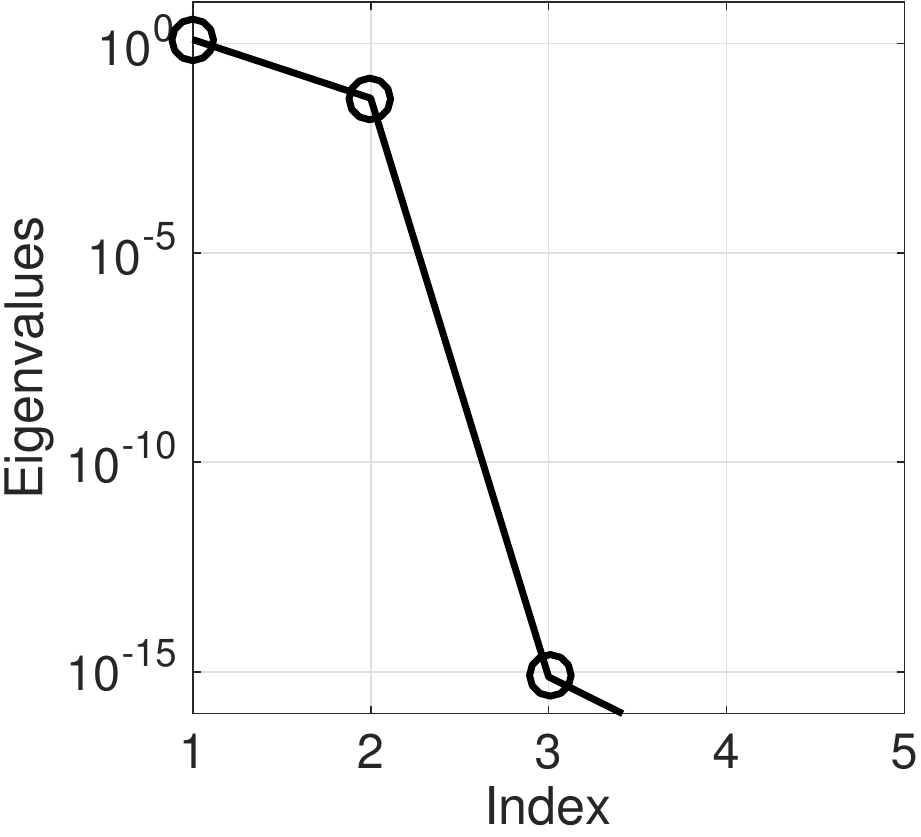}%
}
\hfil
\subfloat[Eigenvectors of $\mC$]{
\label{fig:evecsU}
\includegraphics[width=0.42\textwidth]{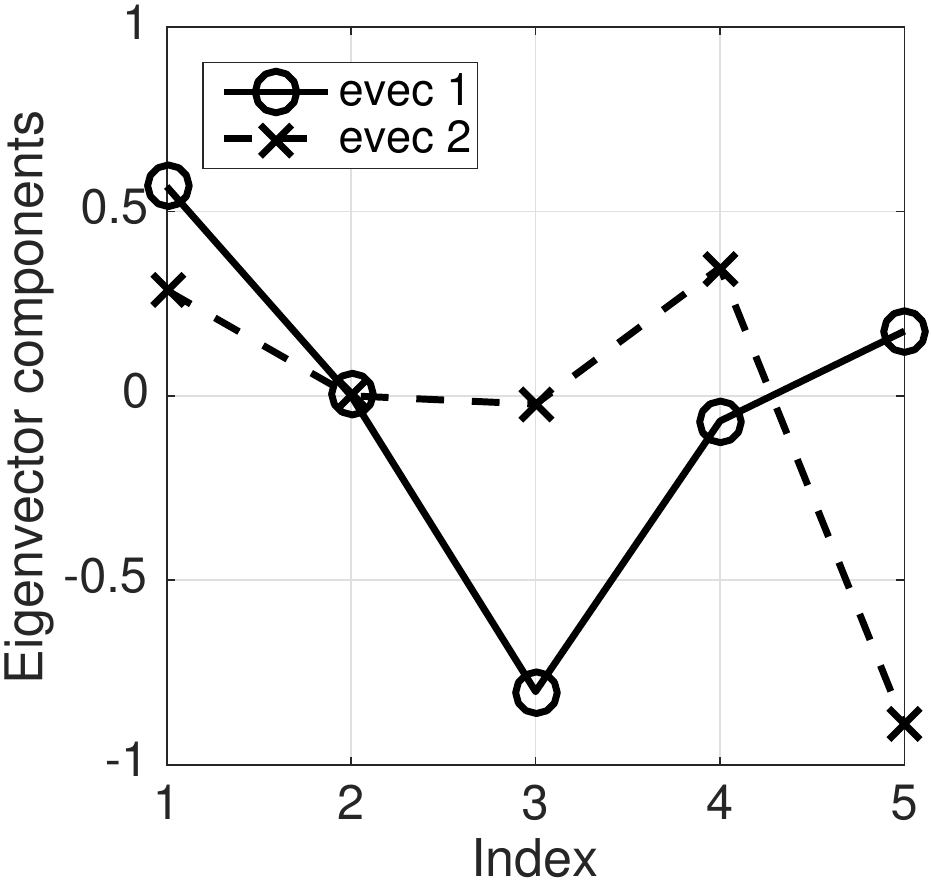}%
}\\
\subfloat[1-D summary plot]{
\label{fig:ssp1U}
\includegraphics[width=0.42\textwidth]{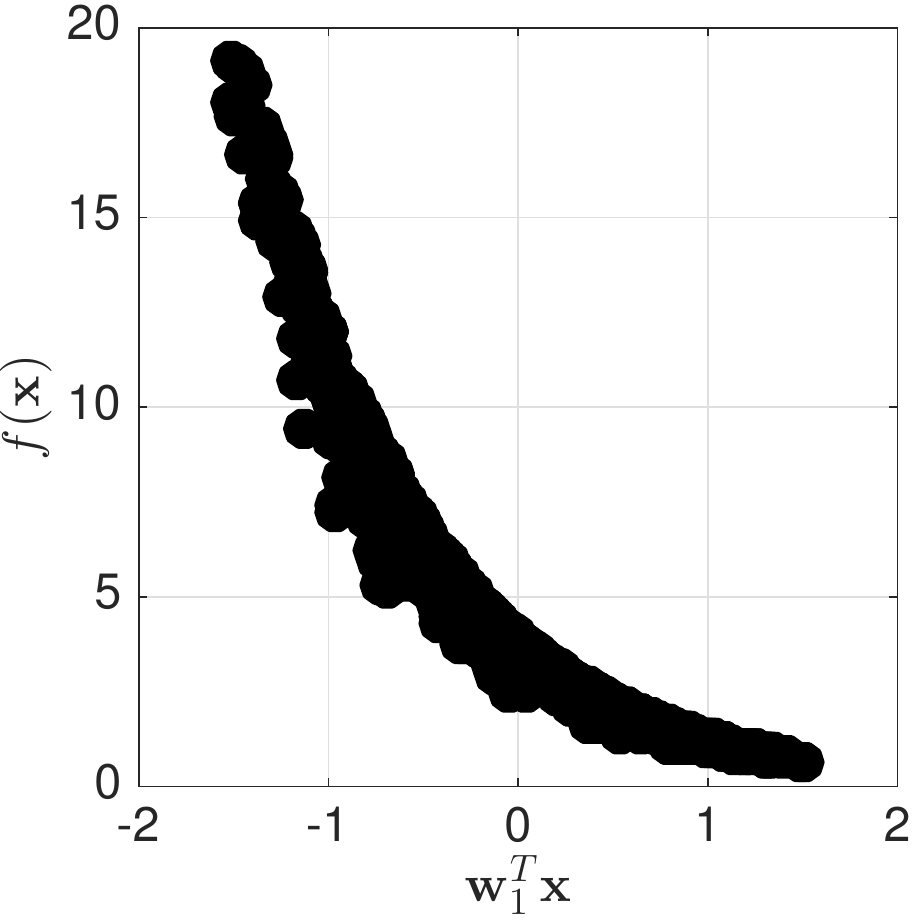}%
}
\hfil
\subfloat[2-D summary plot]{
\label{fig:ssp2U}
\includegraphics[width=0.42\textwidth]{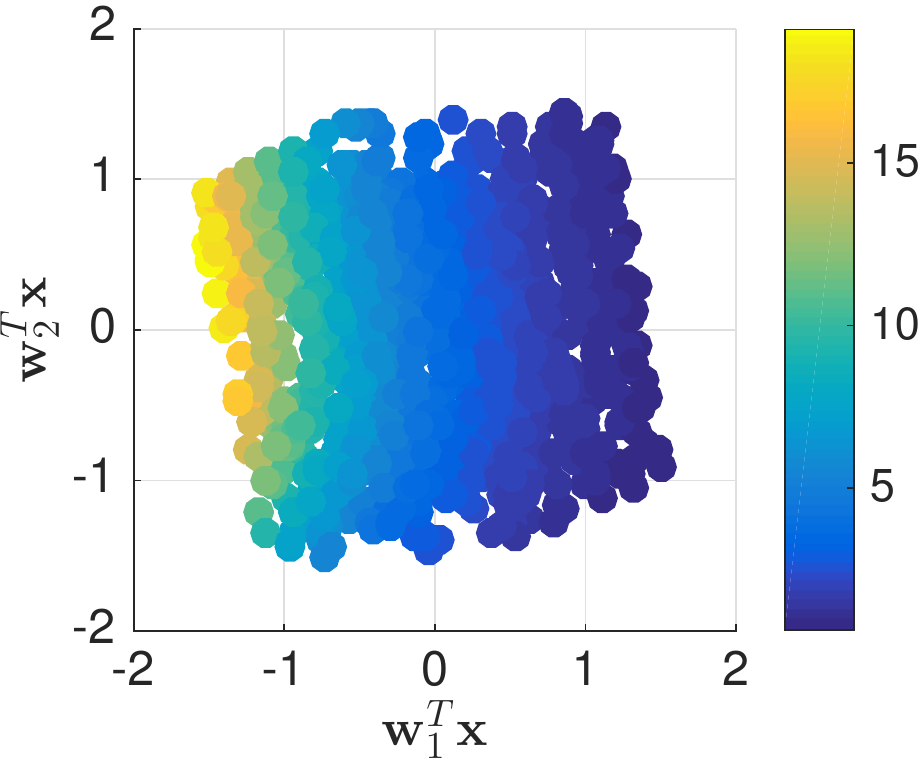}%
}
\caption{These figures represent the active subspace-based dimension reduction for the Hartmann problem's average flow velocity $\uavg$ from \eqref{eq:Uavg_eq_dim}. Figure \ref{fig:evalsU} shows the eigenvalues of $\mC$, and Figure \ref{fig:evecsU} shows the components of $\mC$'s first two eigenvectors. Figures \ref{fig:ssp1U} and \ref{fig:ssp2U} are one- and two-dimensional, respectively, summary plots of the quantity of interest. They reveal the low-dimensional relationship between the two active variables and $\uavg$.}
\label{fig:hartmann-uavg}
\end{figure}

Figure \ref{fig:hartmann-uavg} shows the results of analyzing the active subspace of the Hartmann problem's average flow velocity $\uavg$ from \eqref{eq:Uavg_eq_dim}. Figure \ref{fig:evalsU} shows that all but two eigenvalues are zero (to machine precision). This implies that the active subspace of dimension $n=2$ is sufficient to describe the relationship between the log-transformed inputs and the quantity of interest. Figure \ref{fig:evecsU} shows the components of the first two eigenvectors of $\mC$'s quadrature estimate; the index on the horizontal axis maps to the specific input as in Table \ref{tab:input_ints}. A large eigenvector component reveals that the corresponding parameter is important in defining the active subspace. Notice that both eigenvector components corresponding to $\log(\rho)$ (the second input) are zero; this is consistent with the definition of $\uavg$ in \eqref{eq:Uavg_eq_dim}, which does not depend on fluid density $\rho$. Figure \ref{fig:ssp1U} is a summary plot of 1000 samples of $\uavg$---taken from the quadrature evaluations used to estimate $\mC$---versus the corresponding samples of the active variable. Such plots are commonly used in regression graphics~\citep{cook2009regression}. The plot shows a strong relationship between the first active variable and $\uavg$, so a ridge function of the form \eqref{eq:AS_Red} with one linear combination would be a good approximation; this is validated by the three-order-of-magnitude gap between the first and second eigenvalues. Figure \ref{fig:ssp2U} shows a two-dimensional summary plot with the same data, where the color is the value of $\uavg$, the horizontal axis is the first active variable (defined by the first eigenvector of $\mC$), and the vertical axis is the second active variable (defined by the second eigenvector of $\mC$). Since the eigenvalues with index greater than 2 are zero, the two-dimensional summary plot reveals the complete relationship between the log-transformed inputs and $\uavg$. 

\begin{figure}[!h]
\centering
\subfloat[Eigenvalues of $\mC$]{
\label{fig:evalsB}
\includegraphics[width=0.42\textwidth]{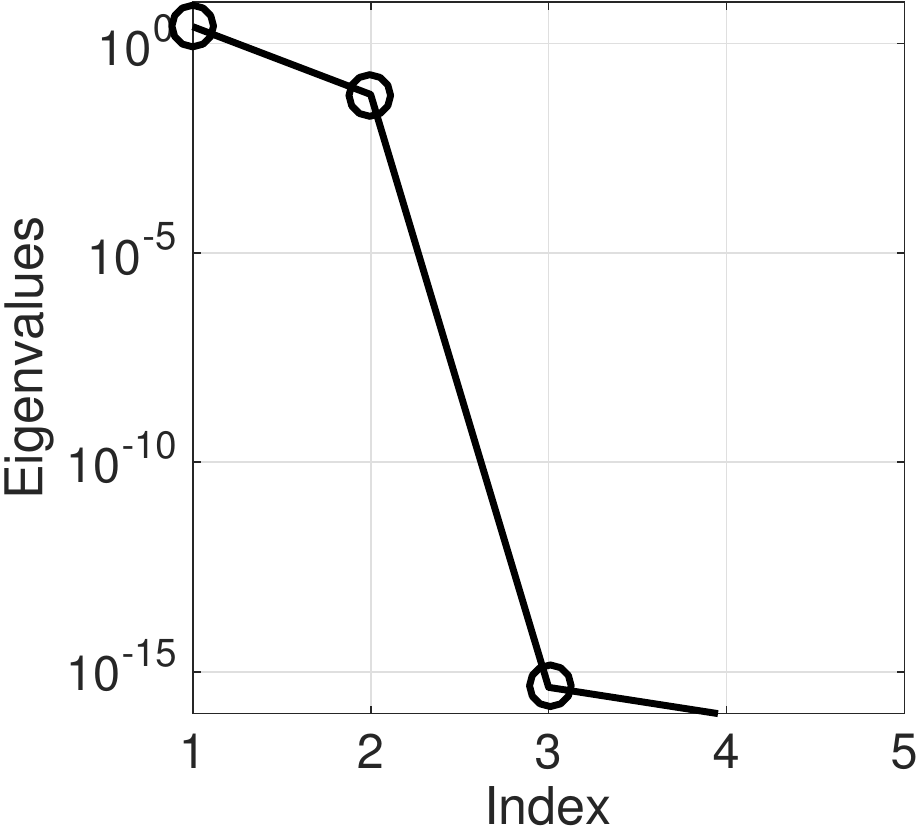}%
}
\hfil
\subfloat[Eigenvectors of $\mC$]{
\label{fig:evecsB}
\includegraphics[width=0.42\textwidth]{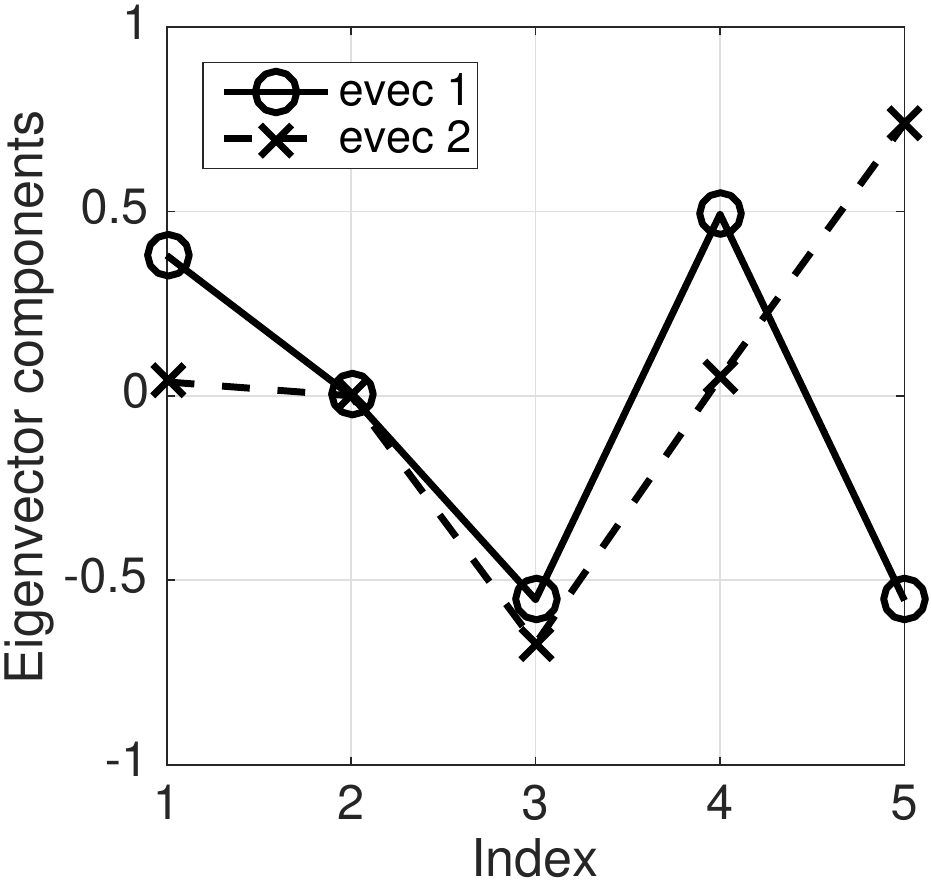}%
}\\
\subfloat[1-D summary plot]{
\label{fig:ssp1B}
\includegraphics[width=0.42\textwidth]{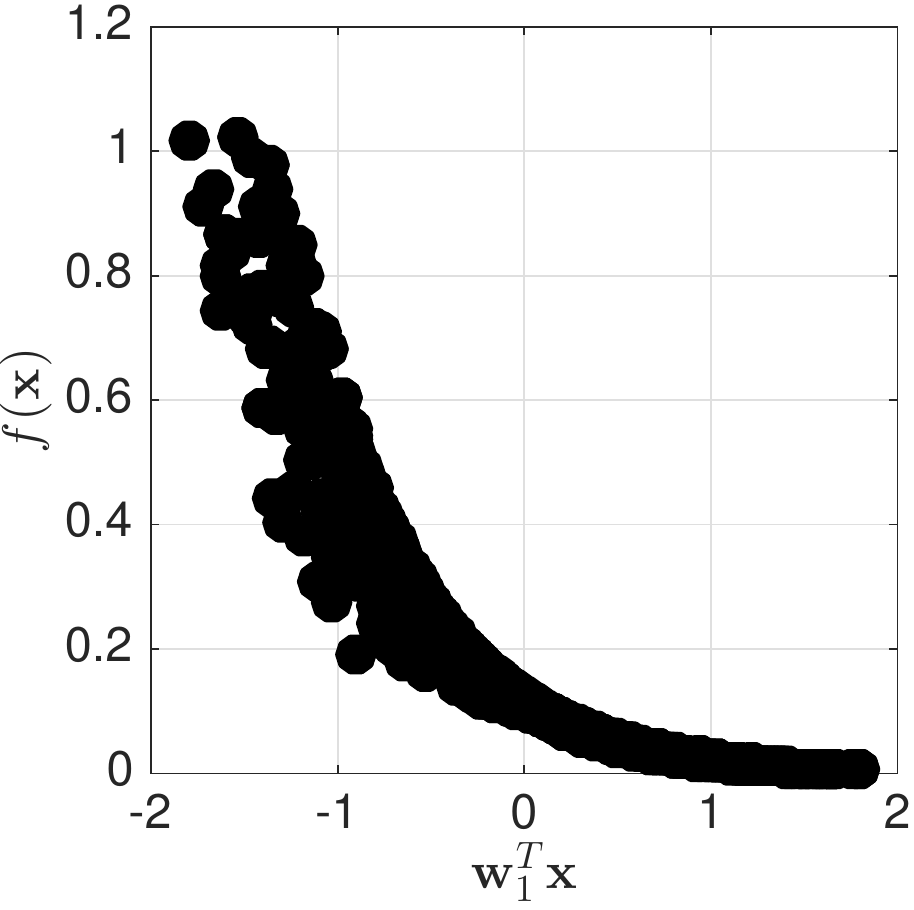}%
}
\hfil
\subfloat[2-D summary plot]{
\label{fig:ssp2B}
\includegraphics[width=0.42\textwidth]{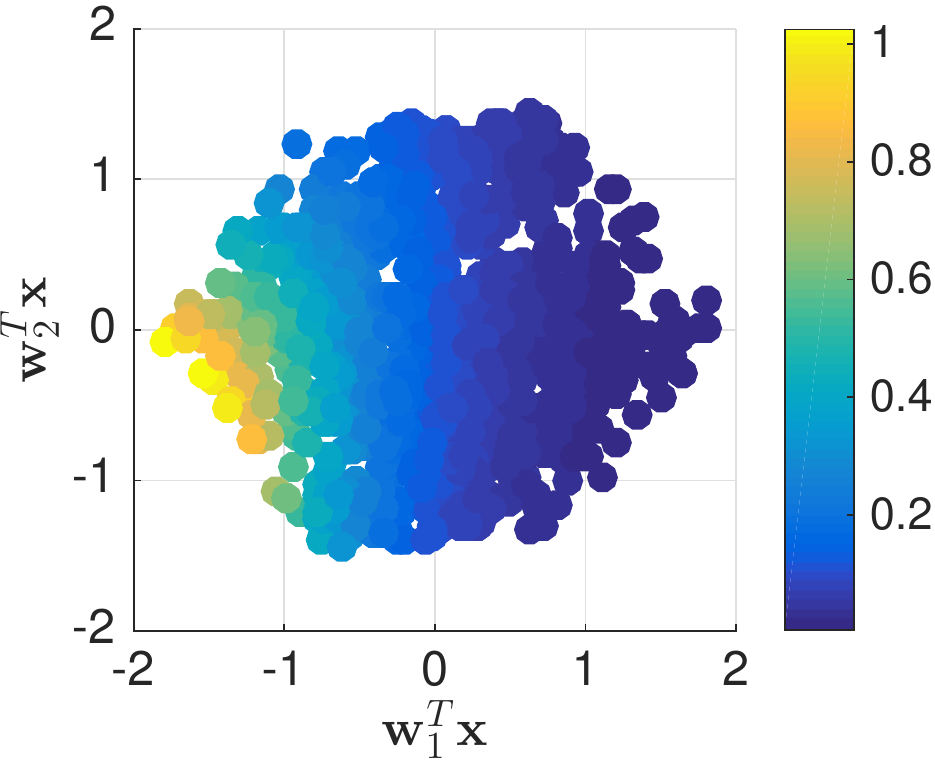}%
}
\caption{These figures represent the active subspace-based dimension reduction for the Hartmann problem's induced magnetic field $\Bind$ from \eqref{eq:Bind_eq_dim}. Figure \ref{fig:evalsB} shows the eigenvalues of $\mC$, and Figure \ref{fig:evecsB} shows the components of $\mC$'s first two eigenvectors. Figures \ref{fig:ssp1B} and \ref{fig:ssp2B} are one- and two-dimensional, respectively, summary plots of the quantity of interest. They reveal the low-dimensional relationship between the two active variables and $\Bind$.}
\label{fig:hartmann-bind}
\end{figure}

Figures \ref{fig:hartmann-bind} show the same information as Figure \ref{fig:hartmann-uavg} for the induced magnetic field $\Bind$ from \eqref{eq:Bind_eq_dim} as the quantity of interest. The plots are similar, and the same comments apply to $\Bind$. 

The dimensional analysis from Section \ref{sec:mhddr} showed that the quantities of interest should depend on 3 unitless quantities. Equation \eqref{eq:nondim_MHD_Eq} expresses the governing equations in terms of the Reynolds number, the Hartmann number, and a dimensionless pressure gradient. With simple scaling, we can write unitless forms of the quantities of interest, $\uavg$ from \eqref{eq:Uavg_eq_dim} and $\Bind$ from \eqref{eq:Bind_eq_dim}, in terms of unitless parameters:
\begin{equation} 
\label{eq:Uavg_eq_dimless}
\uavg^* \;=\; -\frac{\partial p^*_0}{\partial x^*} \frac{\Re}{\Ha^2} \left( 1 - \Ha \, \text{coth} \left( \Ha \right) \right)
\end{equation}
and
\begin{equation} 
\label{eq:Bind_eq_dimless}
\Bind^* \;=\; \frac{\partial p^*_0}{\partial x^*} \frac{\Re}{\Ha} \mu_0^* \left( 1 - \frac{2}{\Ha} \, \text{tanh} \left( \frac{\Ha}{2} \right) \right) .
\end{equation}
Notice the Reynolds number and the dimensionless pressure gradient appear only as a product---which we could define as a new unitless quantity. This explains why the eigenvalues of $\mC$ are zero after the second. 

\subsection{MHD generator problem}
\label{sec:gen}

\begin{figure}[h]
\begin{center}
\includegraphics[width=0.9\textwidth]{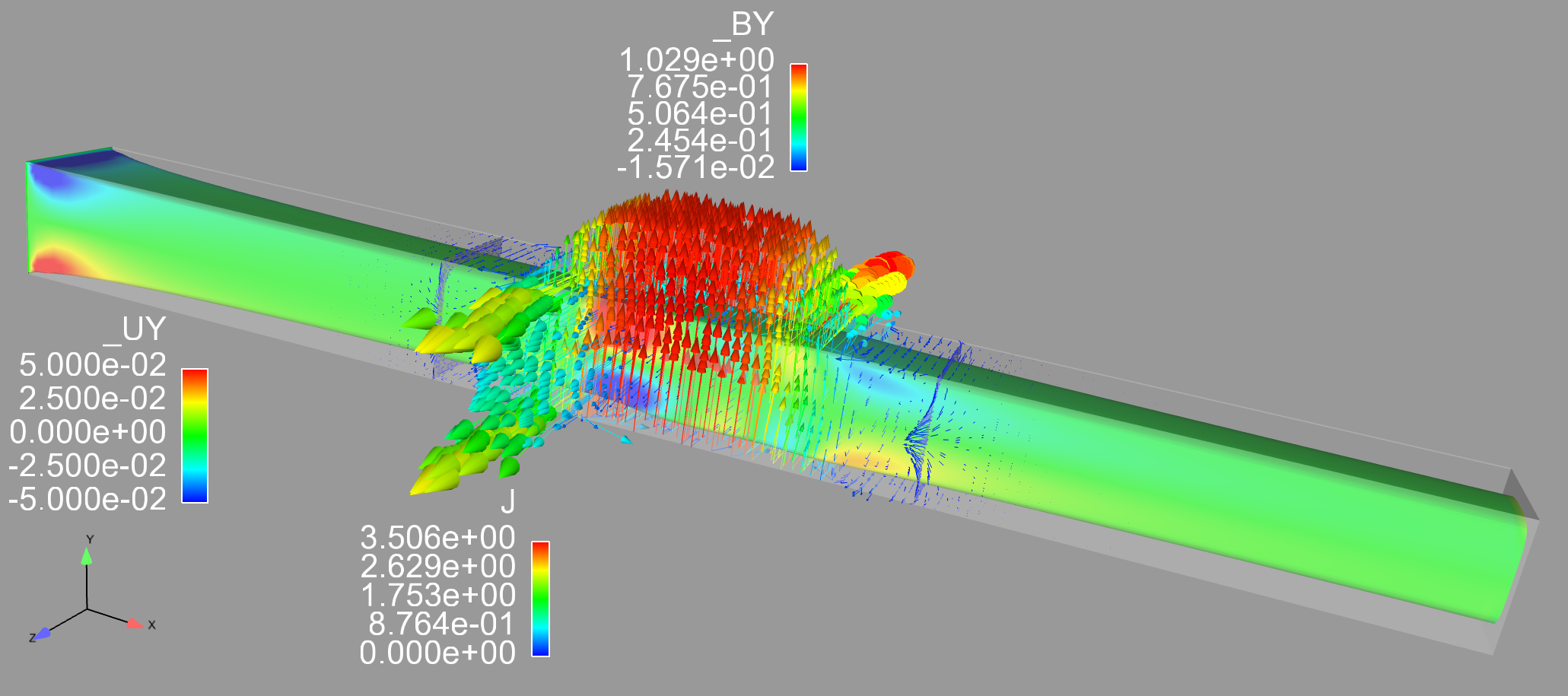}\hspace{0.1in}
\end{center}
\caption{Visualization of flow field from an idealized 3D MHD generator model. The image shows $x$ velocity iso-surface colored by the $y$ velocity. Vectors (colored by magnitude) show vertical magnetic field (applied and induced) and horizontal induced current.}
\label{fig:MHDGeneratorImage}
\end{figure}

This model is a steady-state MHD duct flow configuration representing an idealized MHD generator. The MHD generator induces an electrical current by supplying a set flow-rate of a conducting fluid through an externally supplied vertical magnetic field. The bending of the magnetic field lines produces a horizontal electrical current. The geometric domain for this problem is a square cross-sectional duct of dimensions $15\unit{m}\times 1\unit{m}\times 1\unit{m}$. The simple geometry problem facilitates scalability studies as different mesh sizes can be easily generated. The velocity boundary conditions are set with Dirichlet inlet velocity of $[1,0,0]\,\unit{m}\cdot\unit{s}^{-1}$, no slip on the top, bottom and sides of the channel, and natural boundary conditions on the outflow. The magnetic field boundary conditions on the top and bottom are specified as a set magnetic field configuration $(0,B_y^{gen},0)$, where 
\begin{equation}
B_y^{gen} \;=\;  \frac{1}{2} B_0 \left[ \tanh \left( \frac{x - x_{\text{on}}}{\delta} \right) -  \tanh \left( \frac{x - x_{\text{off}}}{\delta} \right) \right].
\end{equation}
The values $x_{\text{on}}$ and $x_{\text{off}}$ indicate the locations in the $x$ direction where the magnetic field is active. The inlet, outlet, and sides are perfect conductors with ${\bf B} \cdot \hat {\bf n} = {\bf 0}$ and ${\bf E} \times \hat {\bf n} = {\bf 0}$, i.e., the current and magnetic fluxes are zero at these boundaries. Homogeneous Dirichlet boundary conditions are used on all surfaces for the Lagrange multiplier. This problem has similar characteristics to the Hartmann problem with viscous boundary layers and Hartmann layers occurring at the boundaries, and a flow field that is strongly modified by the magnetic field in the section of the duct where it is active. 
Figure \ref{fig:MHDGeneratorImage} shows a solution for this problem for Reynolds number $\Re = 2500$, magnetic Reynolds number $Re_m = 10$, and Hartmann number $\Ha = 5$. The image shows $x$ velocity iso-surface colored by $y$ velocity, where the modification of the inlet constant profile and the parabolic profile at the 
region where the magnetic field is active are evident. Vectors (colored by magnitude) show the vertical magnetic field (applied and induced) and horizontal induced current from the bending of the magnetic field lines.

\begin{table}
\centering
\caption{Indices and intervals for the parameters $\vx$ of the MHD generator problem. These intervals represent the expected operating conditions for the idealized MHD generator.}
\begin{tabular}{llll}
Index & Name & Notation & Interval \\ \hline
1 & fluid viscosity & $\log(\mu)$ & $[\log(0.001), \log(0.01)]$ \\
2 & fluid density & $\log(\rho)$ & $[\log(0.1), \log(10)]$ \\
3 & applied pressure gradient & $\log\left(\frac{\partial p_0}{\partial x}\right)$ & $[\log(0.1), \log(0.5)]$ \\
4 & resistivity & $\log(\eta)$ & $[\log(0.1), \log(10)]$ \\
5 & applied magnetic field & $\log(B_0)$ & $[\log(0.1), \log(1)]$
\end{tabular}
\label{tab:input_ints2}
\end{table}

The fixed physical parameters for the MHD generator are $\mu_0 =1$, $x_{\text{on}} = 4.0$, $x_{\text{off}} = 6.0$, and $\delta = 0.1$. The variable input parameters are the same as in the Hartmann problem. However, the generator uses different input ranges, which can be found in Table \ref{tab:input_ints2}. The probability density function on the space of inputs is a uniform density on the hypercube of log-transformed parameters defined by the ranges in Table \ref{tab:input_ints2}. The quantities of interest are the average flow velocity $\uavg$ and the induced magnetic field $\Bind$, as in the Hartmann problem.

Given values for the input parameters, the MHD generator's solution fields are computed with the Sandia National Laboratory's Drekar multiphysics solver package~\citep{drekar2012}. The package has adjoint capabilities, which enables computation of the derivatives of the quantities of interest with respect to the input parameters~\citep{Shadid2016a}. Each MHD generator model run uses 5.3 CPU-hours (10 minutes on 32 cores), so estimating $\mC$ from \eqref{eq:C} with a tensor product Gauss-Legendre quadrature rule is not possible. Instead, we use a Monte Carlo method to estimate $\mC$ using $M=483$ independent samples from the uniform density on the log-transformed parameters. For details on the accuracy of the Monte Carlo method for estimating active subspaces, see~\cite{constantine2015computing}. 

\begin{figure}[!h]
\centering
\subfloat[Eigenvalues of $\mC$]{
\label{fig:evalsUgen}
\includegraphics[width=0.42\textwidth]{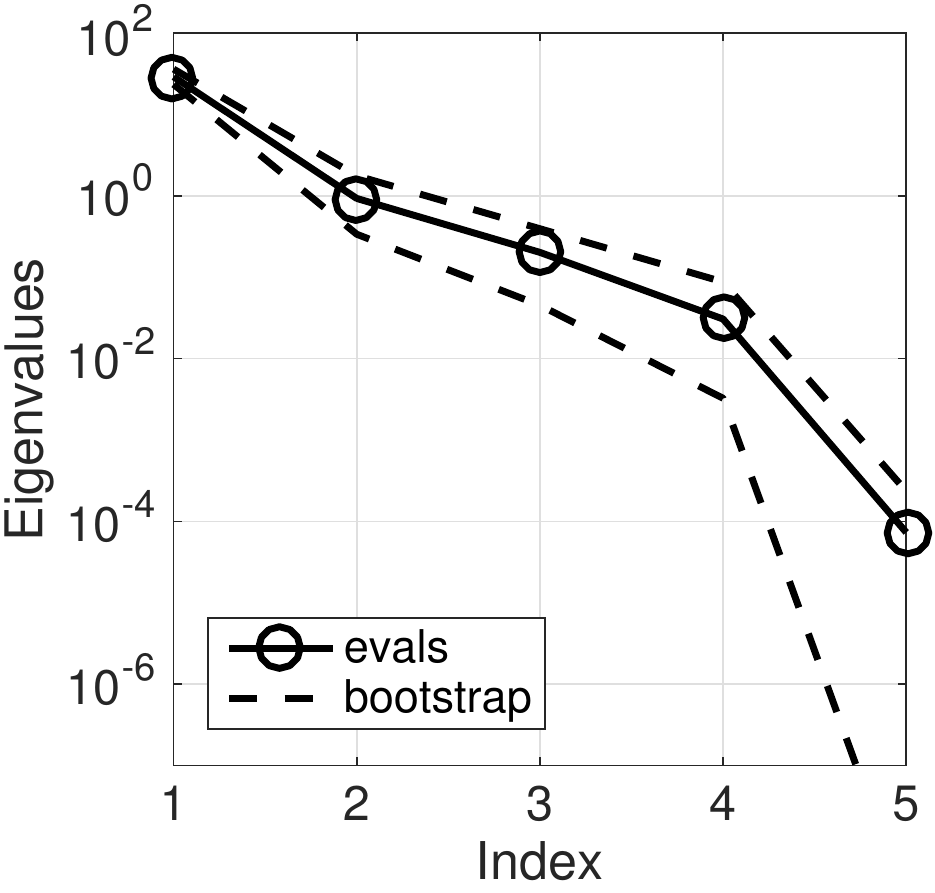}%
}
\hfil
\subfloat[Eigenvectors of $\mC$]{
\label{fig:evecsUgen}
\includegraphics[width=0.42\textwidth]{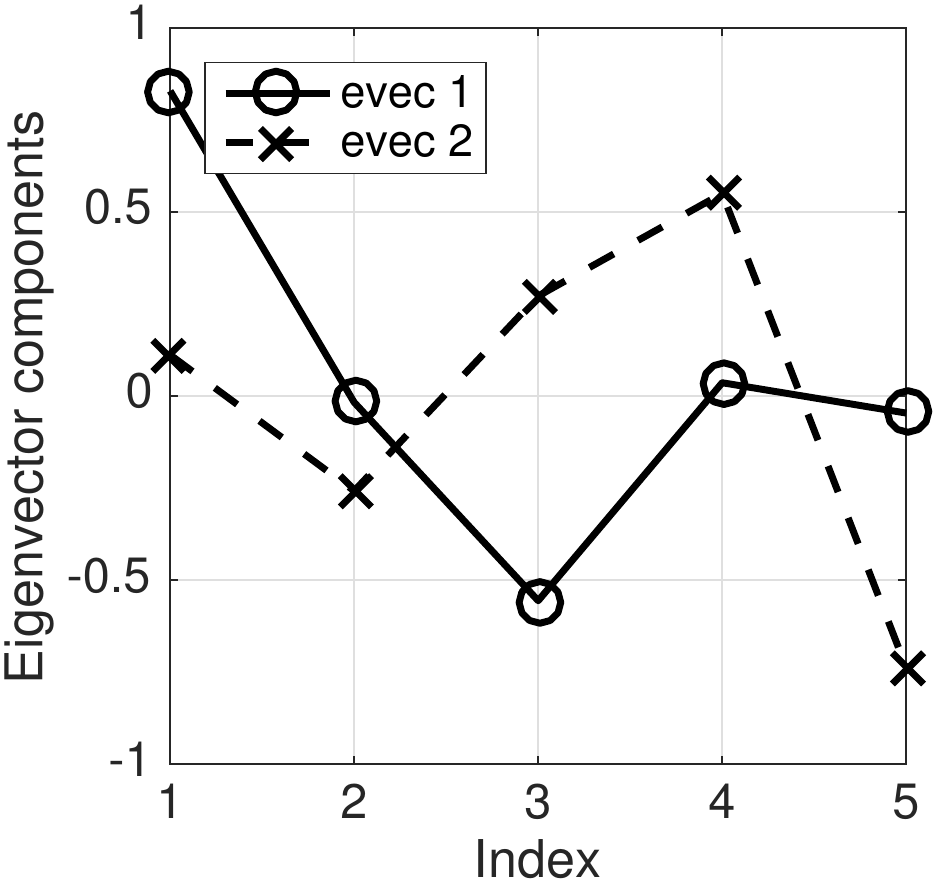}%
}\\
\subfloat[1-D summary plot]{
\label{fig:ssp1Ugen}
\includegraphics[width=0.42\textwidth]{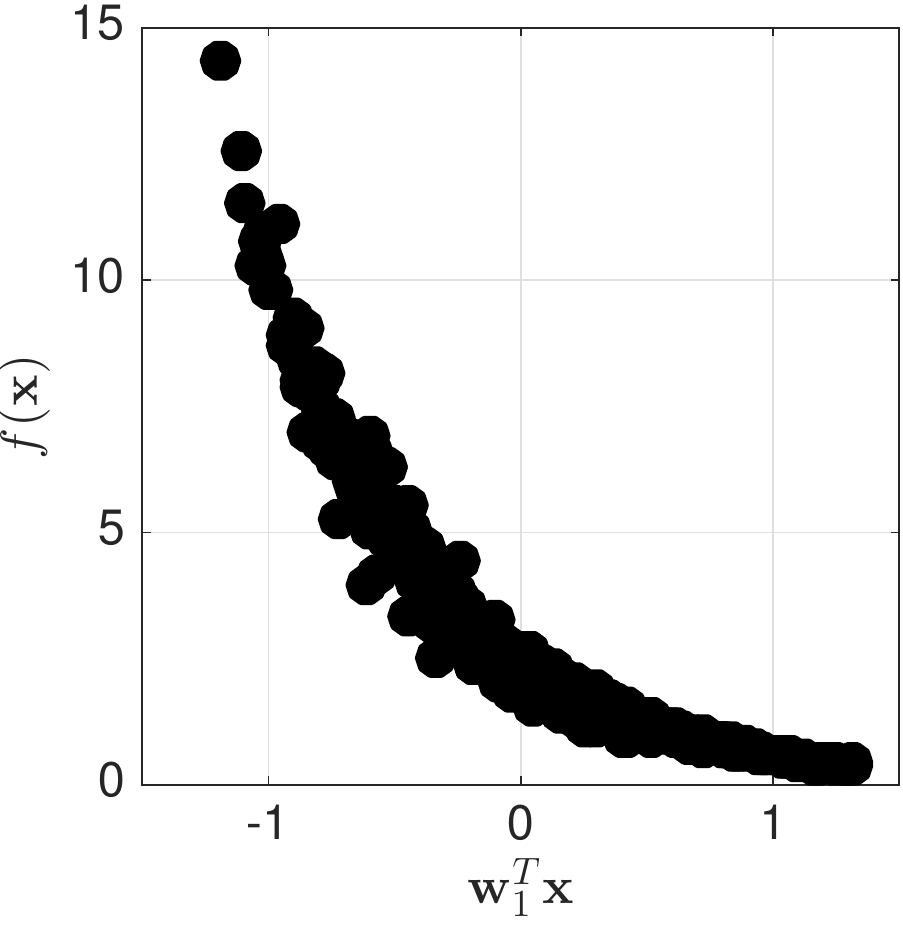}%
}
\hfil
\subfloat[2-D summary plot]{
\label{fig:ssp2Ugen}
\includegraphics[width=0.42\textwidth]{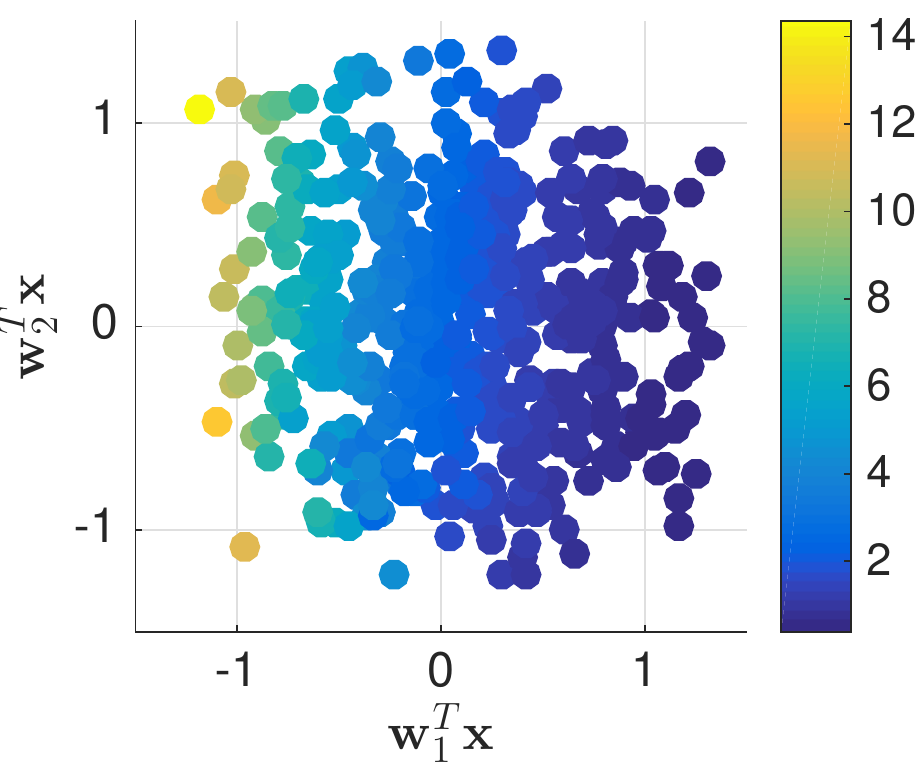}%
}
\caption{These figures represent the active subspace-based dimension reduction for the MHD generator problem's average flow velocity $\uavg$. Figure \ref{fig:evalsUgen} shows the eigenvalues of $\mC$, and Figure \ref{fig:evecsUgen} shows the components of $\mC$'s first two eigenvectors. Figures \ref{fig:ssp1Ugen} and \ref{fig:ssp2Ugen} are one- and two-dimensional, respectively, summary plots of the quantity of interest. They reveal the low-dimensional relationship between the two active variables and $\uavg$.}
\label{fig:gen-uavg}
\end{figure}

Figure \ref{fig:evalsUgen} shows the eigenvalue estimates computed with Monte Carlo for the $\uavg$ quantity of interest. The dashed lines show upper and lower bounds on the eigenvalue estimates computed with a nonparametric bootstrap with 500 bootstrap replicates from the set of 483 gradient samples. We emphasize that since there is no randomness in the map from inputs to outputs (i.e., the computer simulation is deterministic), the bootstrap is a heuristic to estimate the variability due to the Monte Carlo sampling. Estimates of standard error from sample variances are not appropriate, since the eigenvalues are nonlinear functions of the gradient samples. For an example of a similar bootstrap computation for eigenvalues, see~\cite[Chapter 7.2]{efron1994introduction}.  

The fifth eigenvalue is 0.00024\% of the sum of the five eigenvalues, which is consistent with the dimensional analysis from Section \ref{sec:mhddr}---i.e., there should be no more than 4 linear combinations of the model parameters that affect the quantity of interest. And this restriction is reflected in the small fifth eigenvalue. 

\begin{figure}[!h]
\centering
\subfloat[Eigenvalues of $\mC$]{
\label{fig:evalsBgen}
\includegraphics[width=0.42\textwidth]{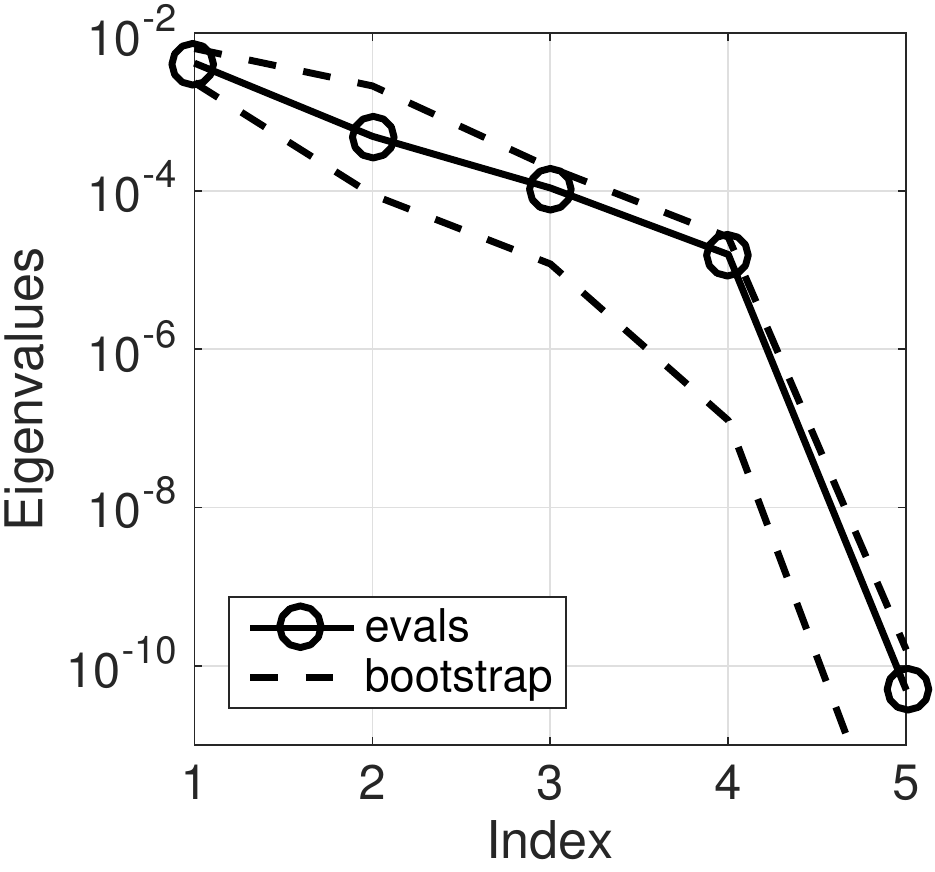}%
}
\hfil
\subfloat[Eigenvectors of $\mC$]{
\label{fig:evecsBgen}
\includegraphics[width=0.42\textwidth]{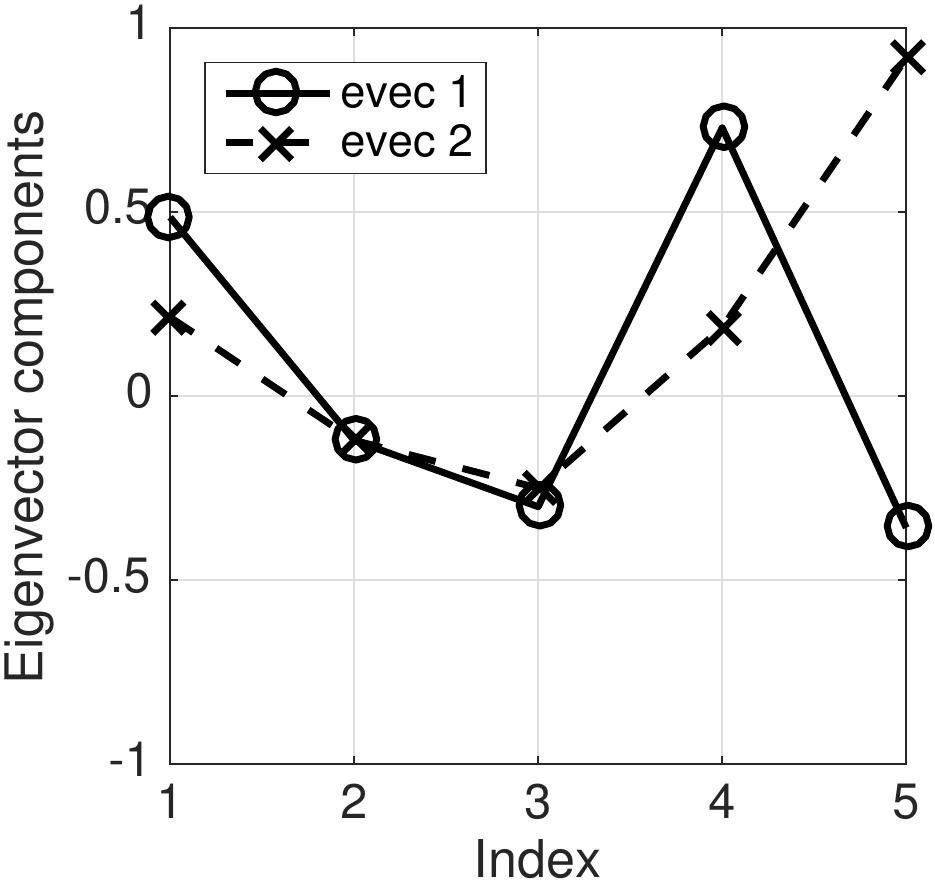}%
}\\
\subfloat[1-D summary plot]{
\label{fig:ssp1Bgen}
\includegraphics[width=0.42\textwidth]{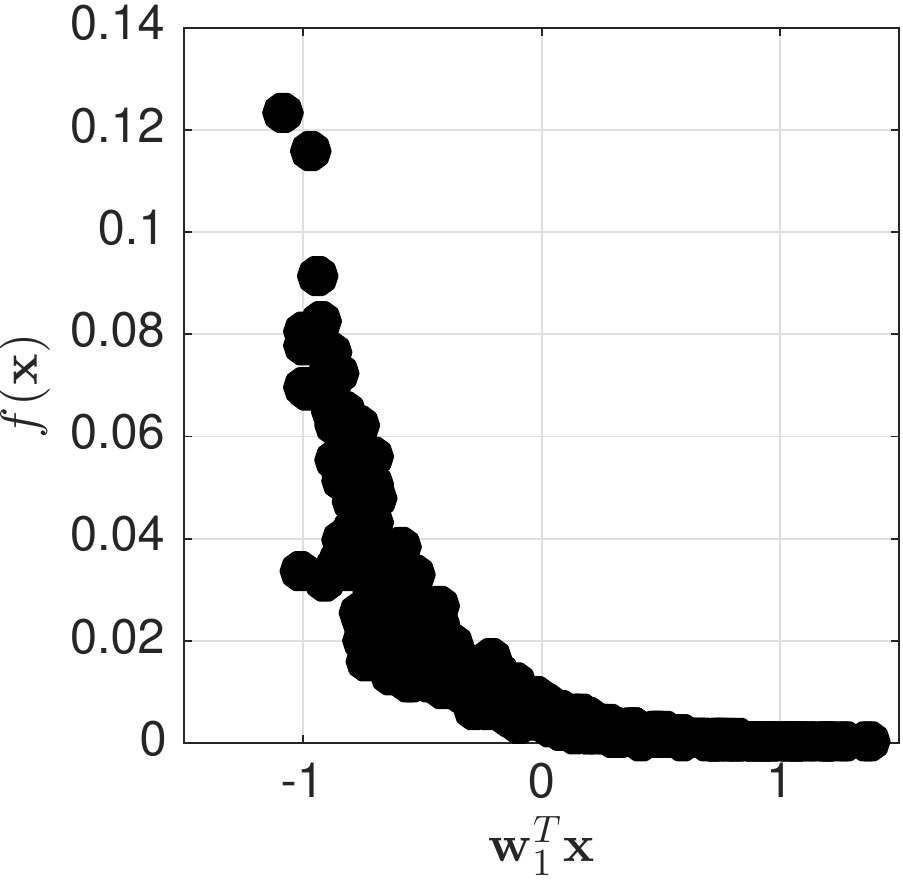}%
}
\hfil
\subfloat[2-D summary plot]{
\label{fig:ssp2Bgen}
\includegraphics[width=0.42\textwidth]{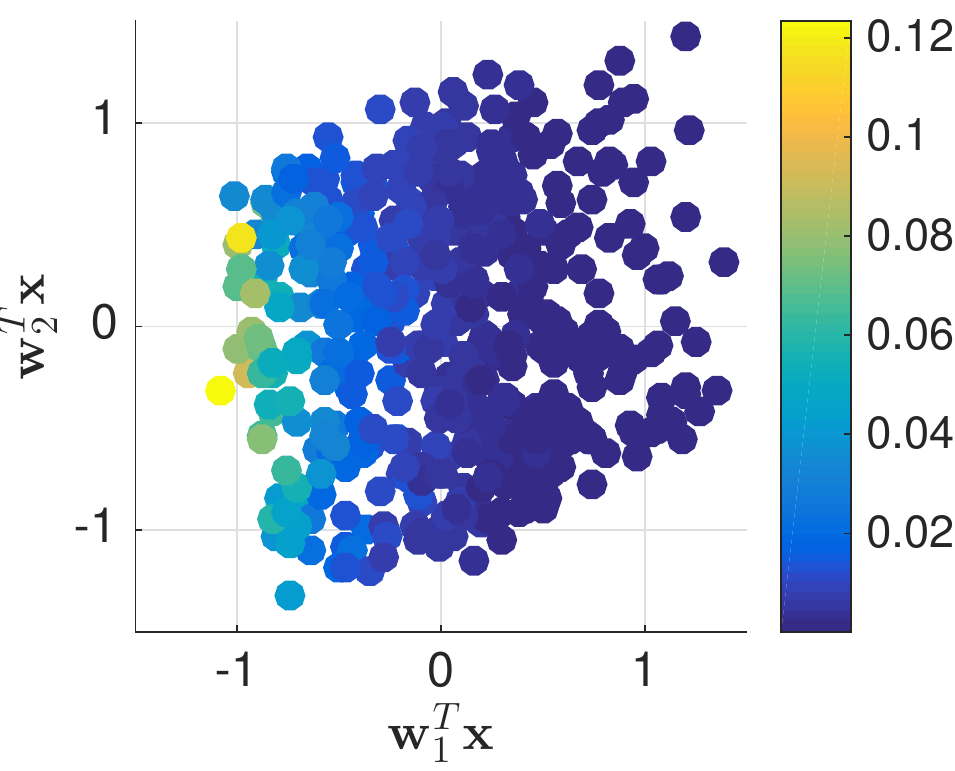}%
}
\caption{These figures represent the active subspace-based dimension reduction for the MHD generator's induced magnetic field $\Bind$. Figure \ref{fig:evalsBgen} shows the eigenvalues of $\mC$, and Figure \ref{fig:evecsBgen} shows the components of $\mC$'s first two eigenvectors. Figures \ref{fig:ssp1Bgen} and \ref{fig:ssp2Bgen} are one- and two-dimensional, respectively, summary plots of the quantity of interest. They reveal the low-dimensional relationship between the two active variables and $\Bind$.}
\label{fig:gen-bind}
\end{figure}

The first two eigenvectors of $\mC$'s Monte Carlo estimate (for $\uavg$) are shown in Figure \ref{fig:evecsBgen}. The magnitudes of the eigenvector components can be used to determine which physical parameters influence the active subspace---i.e., they provide sensitivity information. (See~\cite{diaz2015global} for how to construct sensitivity metrics from active subspaces.) The fluid viscosity $\mu$ and the pressure gradient $\partial p_0/\partial x$ are the most important parameters for the average fluid velocity. This insight agrees with physical intuition, and it is consistent with the same metrics from the Hartmann problem; see Figure \ref{fig:evecsU}. 

Figure \ref{fig:ssp1Ugen} and \ref{fig:ssp2Ugen} show the one- and two-dimensional summary plots for $\uavg$ as a function of the first two active variables using all 483 samples. Similar to Figures \ref{fig:ssp1U} and \ref{fig:ssp2U}, we see a nearly one-dimensional relationship between the log-transformed input parameters and the average velocity, where the one dimension is the first active variable. 

Figure \ref{fig:evalsBgen} shows the eigenvalues for $\mC$'s Monte Carlo estimate, with bootstrap ranges, for the induced magnetic field quantity of interest $\Bind$. In this case, the fifth eigenvalue is 0.000001\% of the sum of the eigenvalues, which is consistent with the dimensional analysis from Section \ref{sec:connections} that shows that any quantity of interest will depend on at most four linear combinations of the log-transformed input parameters. The first eigenvector in Figure \ref{fig:evecsBgen} shows that $\Bind$ depends on all input parameters except the fluid density $\rho$. This is remarkably similar to the dependence seen in the Hartmann problem; see Figure \ref{fig:evecsB}. 

The one- and two-dimensional summary plots for $\Bind$ are in Figures \ref{fig:ssp1Bgen} and \ref{fig:ssp2Bgen}. There appears to be a region in the parameter space---when the first active variable is positive---where the relationship between the inputs and $\Bind$ is well characterized by one linear combination of the log-transformed inputs. However, the one-dimensional character of that relationship degrades as the first active variable decreases. (Note that the first eigenvector is only unique up to a sign, so this relationship could be inverted. To connect to the model's input parameters, one must examine the signs of the individual eigenvector components. For example, the fourth and fifth components of the first eigenvector have opposite signs. Thus, the corresponding model parameters---resistivity $\eta$ and applied magnetic field $B_0$---affect $\Bind$ in opposite directions, on average.) 

\section{Conclusions}

We have reviewed two methods for dimension reduction in physical systems: (i) dimensional analysis that uses the physical quantities' units and (ii) active subspaces that use the gradient of the output with respect to the inputs to identify an important low-dimensional subspace of the input space. We also reviewed the connection between these two methods via a log transform of the input parameters---namely, that the dimensional analysis provides an upper bound on the number of linear combinations of log-transformed parameters that control any quantity of interest. 

We applied these techniques to two quantities of interest---average flow velocity and induced magnetic field---from two magnetohydrodynamics models that apply to power generation: (i) the Hartmann problem that admits closed form expressions for the quantities of interest and (ii) a large-scale computational model of coupled fluid flow, magnetic fields, and electric current in a three-dimensional duct. The computational model has adjoint capabilities that enable gradient evaluations. 

The insights from the active subspace are consistent with the dimensional analysis. In particular, there are at most four linear combinations of log-transformed parameters that affect the quantity of interest---which is a reduction from an ambient dimension of 5 to an intrinsic dimension of 4. The Hartmann problem has a further reduction to an intrinsic dimension of 2---i.e., two linear combinations of log-transformed parameters are sufficient to characterize the quantities of interest. The eigenvalues of the matrix $\mC$ that defines the active subspace rank the importance of the linear combinations, which offers more insight into the input/output relationships than the dimensional analysis alone. 

The eigenvector components that define the active subspace reveal the sensitivities of the quantities of interest with respect to the input parameters. These metrics are consistent between the Hartmann problem and the large-scale generator model. In particular, (i) the average flow velocity depends mainly on the fluid viscosity and the applied pressure gradient, and (ii) the induced magnetic field depends strongly on all parameters except the fluid density. 

The summary plots with the first and second eigenvectors show the low-dimensional relationship between the log-transformed inputs and the quantities of interest. In particular, the average flow velocity can be well approximated by a function of one linear combination of the log-transformed inputs in both the Hartmann problem and large-scale generator. The same is true for the induced magnetic field from the Hartmann problem. However, the induced magnetic field in the large-scale generator admits a one-dimensional characterization in a subset of the parameter domain. By combining the information in the one-dimensional summary plot with the sensitivities reveals by the first eigenvector's components, we conclude that in regions of (i) low fluid viscosity, (ii) high applied pressure gradient, (iii) low resistivity, and (iv) high applied magnetic field, the relationship between the log-transformed inputs and the induced magnetic field requires more than two linear combinations of input parameters for accurate approximation. 

These insights offer directions for further development in MHD models for effective power generation. Additionally, the revealed relationships between the inputs and outputs will assist in (i) quantifying parametric uncertainties and (ii) enabling computational design for MHD power generation.



\section*{Acknowledgments}

The first author's work is supported by the Ben L.~Fryrear Ph.D.~Fellowship in Computational Science. The second author's work supported by the U.S. Department of Energy Office of Science, Office of Advanced Scientific Computing Research, Applied Mathematics program under Award Number DE-SC-0011077. Sandia National Laboratories is a multi-program laboratory managed and operated by Sandia Corporation, a wholly owned subsidiary of Lockheed Martin Corporation, for the U.S. Department of Energy’s National Nuclear Security Administration under contract DE-AC04-94AL85000.

\bibliographystyle{natbib}
\bibliography{mhd_as}

\begin{thebibliography}{}

\bibitem[Barenblatt(1996)Barenblatt]{Barenblatt1996}
Barenblatt, G.~I. (1996).
\newblock {\em Scaling, Self-similarity, and Intermediate Asymptotics\/}.
\newblock Cambridge University Press, Cambridge.

\bibitem[Constantine and Gleich(2015)Constantine and
  Gleich]{constantine2015computing}
Constantine, P. and Gleich, D. (2015).
\newblock Computing active subspaces with {Monte Carlo}.
\newblock {\em arXiv:1408.0545v2\/}.

\bibitem[Constantine {\em et~al.}(2014)Constantine, Dow, and
  Wang]{constantine2014active}
Constantine, P., Dow, E., and Wang, Q. (2014).
\newblock Active subspace methods in theory and practice: Applications to
  kriging surfaces.
\newblock {\em SIAM Journal on Scientific Computing\/}, {\bf 36}(4),
  A1500--A1524.

\bibitem[Constantine {\em et~al.}(2015a)Constantine, Emory, Larsson, and
  Iaccarino]{Constantine2015}
Constantine, P., Emory, M., Larsson, J., and Iaccarino, G. (2015a).
\newblock Exploiting active subspaces to quantify uncertainty in the numerical
  simulation of the {HyShot II} scramjet.
\newblock {\em Journal of Computational Physics\/}, {\bf 302}, 1--20.

\bibitem[Constantine {\em et~al.}(2016)Constantine, del Rosario, and
  Iaccarino]{constantine2016many}
Constantine, P., del Rosario, Z., and Iaccarino, G. (2016).
\newblock Many physical laws are ridge functions.
\newblock {\em arXiv:1605.07974\/}.

\bibitem[Constantine(2015)Constantine]{asmbook}
Constantine, P.~G. (2015).
\newblock {\em Active Subspaces: Emerging Ideas for Dimension Reduction in
  Parameter Studies\/}.
\newblock SIAM, Philadelphia.

\bibitem[Constantine {\em et~al.}(2015b)Constantine, Zaharatos, and
  Campanelli]{constantine2015discovering}
Constantine, P.~G., Zaharatos, B., and Campanelli, M. (2015b).
\newblock Discovering an active subspace in a single-diode solar cell model.
\newblock {\em Statistical Analysis and Data Mining: The ASA Data Science
  Journal\/}, {\bf 8}(5-6), 264--273.

\bibitem[Cook(2009)Cook]{cook2009regression}
Cook, R.~D. (2009).
\newblock {\em Regression graphics: Ideas for studying regressions through
  graphics\/}.
\newblock John Wiley \& Sons, 2nd edition.

\bibitem[{Cowling} and {Lindsay}(1957){Cowling} and {Lindsay}]{Cowl:57}
{Cowling}, T.~G. and {Lindsay}, R.~B. (1957).
\newblock Magnetohydrodynamics.
\newblock {\em Physics Today\/}, {\bf 10}, 40.

\bibitem[Diaz and Constantine(2015)Diaz and Constantine]{diaz2015global}
Diaz, P. and Constantine, P.~G. (2015).
\newblock Global sensitivity metrics from active subspaces.
\newblock {\em arXiv:1510.04361v1\/}.

\bibitem[Efron and Tibshirani(1994)Efron and Tibshirani]{efron1994introduction}
Efron, B. and Tibshirani, R.~J. (1994).
\newblock {\em An Introduction to the Bootstrap\/}.
\newblock CRC press.

\bibitem[Jefferson {\em et~al.}(2015)Jefferson, Gilbert, Constantine, and
  Maxwell]{Jefferson2015}
Jefferson, J.~L., Gilbert, J.~M., Constantine, P.~G., and Maxwell, R.~M.
  (2015).
\newblock Active subspaces for sensitivity analysis and dimension reduction of
  an integrated hydrologic model.
\newblock {\em Computers \& Geosciences\/}, {\bf 83}, 127--138.

\bibitem[Lin {\em et~al.}(2010)Lin, Shadid, Tuminaro, Sala, Hennigan, and
  Pawlowski]{Lin2010}
Lin, P.~T., Shadid, J.~N., Tuminaro, R.~S., Sala, M., Hennigan, G.~L., and
  Pawlowski, R.~P. (2010).
\newblock A parallel fully coupled algebraic multilevel preconditioner applied
  to multiphysics {PDE} applications: Drift-diffusion, flow/transport/reaction,
  resistive {MHD}.
\newblock {\em International Journal for Numerical Methods in Fluids\/}, {\bf
  64}(10-12), 1148--1179.

\bibitem[{National Institute of Standards and Technology}(2016){National
  Institute of Standards and Technology}]{nistunits}
{National Institute of Standards and Technology} (2016).
\newblock International system of units (si).
\newblock \url{http://physics.nist.gov/cuu/Units/units.html}.
\newblock Accessed: 2016-06-11.

\bibitem[NETL(2014)NETL]{mhdworkshop}
NETL (2014).
\newblock Magnetohydrodynamic power generation workshop summary report.
\newblock Technical report, National Energy Technology Laboratory.

\bibitem[Pawlowski {\em et~al.}(2012)Pawlowski, Shadid, Smith, Cyr, and
  Weber]{drekar2012}
Pawlowski, R., Shadid, J., Smith, T., Cyr, E., and Weber, P. (2012).
\newblock Drekar::cfd-a turbulent fluid-flow and conjugate heat transfer code:
  Theory manual version 1.0.
\newblock Technical Report SAND2012-2697, Sandia National Laboratories.

\bibitem[Pinkus(2015)Pinkus]{pinkus2015}
Pinkus, A. (2015).
\newblock {\em Ridge Functions\/}.
\newblock Cambridge University Press, Cambridge.

\bibitem[Shadid {\em et~al.}(2010)Shadid, Pawlowski, Banks, Chacón, Lin, and
  Tuminaro]{Shadid2010}
Shadid, J., Pawlowski, R., Banks, J., Chacón, L., Lin, P., and Tuminaro, R.
  (2010).
\newblock Towards a scalable fully-implicit fully-coupled resistive {MHD}
  formulation with stabilized {FE} methods.
\newblock {\em Journal of Computational Physics\/}, {\bf 229}(20), 7649--7671.

\bibitem[Shadid {\em et~al.}(2016a)Shadid, Pawlowski, Cyr, Tuminaro,
  Chac\'{o}n, and Weber]{Shadid2016}
Shadid, J., Pawlowski, R., Cyr, E., Tuminaro, R., Chac\'{o}n, L., and Weber, P.
  (2016a).
\newblock Scalable implicit incompressible resistive {MHD} with stabilized {FE}
  and fully-coupled {Newton–Krylov-AMG}.
\newblock {\em Computer Methods in Applied Mechanics and Engineering\/}, {\bf
  304}, 1--25.

\bibitem[Shadid {\em et~al.}(2016b)Shadid, Smith, Cyr, Wildey, and
  Pawlowski]{Shadid2016a}
Shadid, J., Smith, T., Cyr, E., Wildey, T., and Pawlowski, R. (2016b).
\newblock Stabilized {FE} simulation of prototype thermal-hydraulics problems
  with integrated adjoint-based capabilities.
\newblock {\em Journal of Computational Physics\/}.
\newblock In press.

\bibitem[Sondak {\em et~al.}(2015)Sondak, Shadid, Oberai, Pawlowski, Cyr, and
  Smith]{Sondak}
Sondak, D., Shadid, J., Oberai, A., Pawlowski, R., Cyr, E., and Smith, T.
  (2015).
\newblock A new class of finite element variational multiscale turbulence
  models for incompressible magnetohydrodynamics.
\newblock {\em Journal of Computational Physics\/}, {\bf 295}, 596--616.

\bibitem[USGAO(1993)USGAO]{gao93}
USGAO (1993).
\newblock {FOSSIL FUELS: The Department of Energy's Magnetohydrodynamics
  Development Program}.
\newblock Technical Report RCED-93-174, U.S. Government Accountability Office.

\end{thebibliography}

\end{document}